\documentclass[nonblindrev]{informs3}

\OneAndAHalfSpacedXI

\usepackage{natbib}
 \bibpunct[, ]{(}{)}{,}{a}{}{,}%
 \def\bibfont{\small}%
\usepackage{tikz}
\usepackage{pgfplots}
\usepackage{graphicx}
\usepackage{subcaption}
\graphicspath{ {./images/} }
\usepackage{amsmath}
\usepackage{multirow}
\usepackage{diagbox}
\usepackage{float}

\usepackage{multirow}
\usepackage{algorithm}
\usepackage[noend]{algorithmic}

\newtheorem{example}{Example}
\newtheorem{lem}{Lemma}
\newtheorem{prop}{Proposition}
\newtheorem{insight}{Insight}
\newtheorem{theorem}{Theorem}

\newtheorem{remark}{Remark}
\newtheorem{obs}{Observation}
\usepackage{array}

\usepackage{tikz,pgfplots}
\usetikzlibrary{matrix}
\usepgfplotslibrary{groupplots}
\pgfplotsset{compat=newest}
\usetikzlibrary{calc, positioning, shapes, fit, arrows, shadows, chains}
\tikzstyle{line} = [ draw, -latex']
\usetikzlibrary{shapes, arrows}
\usepgfplotslibrary{units}
\usepackage{url}

\usepackage{booktabs}

\usepackage{color}

\definecolor{color1}{rgb}{1.0, 0.5, 0.0}
\definecolor{color2}{rgb}{0.9, 0.4, 0.3}
\definecolor{color3}{rgb}{0.8, 0.3, 0.5}
\definecolor{color4}{rgb}{0.6, 0.2, 0.8}
\definecolor{color5}{rgb}{0.4, 0.1, 0.9}
\definecolor{color6}{rgb}{0.2, 0.0, 1.0}

\begin{document}

\RUNAUTHOR{Kayac{\i}k et al.}
\RUNTITLE{Adaptive Multistage Stochastic Programming}

\TITLE{Adaptive Multistage Stochastic Programming}



\ARTICLEAUTHORS{%
	\AUTHOR{Sezen Ece Kayac{\i}k}
	\AFF{Department of Operations, Faculty of Economics and Business, University of Groningen, the Netherlands, \EMAIL{s.e.kayacik@rug.nl}}
 \AUTHOR{Beste Basciftci}
	\AFF{Department of Business Analytics, Tippie College of Business, University of Iowa, Iowa City, Iowa, the United States, \EMAIL{beste-basciftci@uiowa.edu}}
 \AUTHOR{Albert H. Schrotenboer}
	\AFF{Operations, Planning, Accounting \& Control Group, School of Industrial Engineering, Eindhoven University of Technology, the Netherlands,  \EMAIL{a.h.schrotenboer@tue.nl}}
 \AUTHOR{Evrim Ursavas}
	\AFF{Department of Operations, Faculty of Economics and Business, University of Groningen, the Netherlands, \EMAIL{e.ursavas@rug.nl}}
 }
\ABSTRACT{%
Multistage stochastic programming is a powerful tool allowing decision-makers to revise their decisions at each stage based on the realized uncertainty. However, in practice, organizations are not able to be fully flexible, as decisions cannot be revised too frequently due to their high organizational impact. Consequently, decision commitment becomes crucial to ensure that initially made decisions remain unchanged for a certain period. This paper introduces adaptive multistage stochastic programming, a new optimization paradigm that strikes an optimal balance between decision flexibility and commitment by determining the best stages to revise decisions depending on the allowed level of flexibility. We introduce a novel mathematical formulation and theoretical properties eliminating certain constraint sets. Furthermore, we develop a decomposition method that effectively handles mixed-integer adaptive multistage programs by adapting the integer L-shaped method and Benders decomposition. Computational experiments on stochastic lot-sizing and generation expansion planning problems show substantial advantages attained through optimal selections of revision times when flexibility is limited, while demonstrating computational efficiency of the proposed properties and solution methodology. Optimizing revision times in a less flexible case can outperform arbitrary selection in a more flexible case. By adhering to these optimal revision times, organizations can achieve performance levels comparable to fully flexible settings.

}



\maketitle

\vspace{-10mm}
\section{Introduction}
Organizations frequently need to make decisions with lasting impact while facing an uncertain future. This is complex, in particular, when such decisions are made sequentially because every decision taken `now' should anticipate the impact of future information and decisions \citep{powell2019unified}. Theoretically, high-quality decisions require \textit{decision flexibility}: the ability to revise earlier made decisions to act upon the latest information, which within a business context translates into `agility' \citep{hbr2022}. However, many organizations cannot adopt agile practices because many decisions cannot be revised too frequently due to their high organizational impact. This necessitates \textit{decision commitment} to assure that initially made decisions remain unchanged for a certain amount of time despite new information or evolving circumstances. In practice, decision-makers must balance decision flexibility with decision commitment; they limit the number of revisions and identify \textit{when} to revise in advance of the actual revision decision. The right balance will result in revising decisions only at the most urgent times, while otherwise committing to decisions previously made. In this study, we provide the first general stochastic optimization approach that determines the optimal revision times to balance flexibility and commitment in stochastic decision-making.

Decision flexibility plays a significant role in effectively addressing challenges in today's uncertain business environment \citep{hbr2010,2hbr2022}. However, it is crucial to acknowledge that frequently revising an earlier-made decision (based on new information) is not always favorable or feasible in real life. This is acknowledged in supply chain management \citep{wadhwa2008framework}, economics \citep{rudloff2014time}, and policy-making \citep{chemama2019consumer}. For example, in supply chain management, producers often have constraints that limit frequently adjusting production quantities and schedules. These limitations can arise from various factors such as production capacity, high setup costs, the need to allocate sources among multiple customers, and predetermined production schedules for other products that cannot be easily changed \citep{kouvelis2006supply, milner2005order}. Besides, a recent survey by \cite{McKinsey2019} highlights the time-consuming nature of decision-making, revealing that respondents spend over 30\% of their working hours on this task, much of which is used ineffectively. For such time-consuming decision revisions, businesses may opt for limiting the number of revisions. In other words, finding the right balance between decision flexibility and commitment is essential. 

As another relevant area, the challenges faced in the energy sector, specifically through the transition to climate-neutral energy systems, necessitate a balance between decision flexibility and commitment. Long-term energy policy decisions must effectively incorporate the inherent uncertainties of the energy system, technological developments, and societal changes \citep{kohler2019agenda}. \cite{meadowcroft2009politics} emphasizes the importance of flexibility and revising decisions based on new information while, at the same time,  ensuring reliability to investors and stakeholders during long-term energy transitions. Accordingly, policymakers should announce their policies and potential revision timelines years in advance. This, for example, is required in developing long-term transition plans toward a green hydrogen economy in the Northern Netherlands, where Europe's first hydrogen valley is being built \citep{newenergycoalition2020}. In such strategic plans, which involve costly investment decisions and multiple stakeholders, frequent changes are undesirable due to the potential disruptions they may cause. Additionally, contractual limitations or the lengthy process required to establish the necessary infrastructure for the energy system can limit flexibility \citep{munoz2013engineering,careri2011generation}, making commitments requisite during planning. \cite{latorre2003classification} emphasize the need for theory and tools that consider the limited flexibility level of energy expansion planning process to adapt to uncertainties in energy markets.

While ignoring the need for decision commitment, many optimization paradigms allow for \textit{full} decision flexibility. Multistage stochastic programming is, arguably, one of most popular approaches amidst such optimization paradigms \citep{birge1985decomposition, daryalal2022lagrangian, bertsimas2023data}; it facilitates sequential decision-making under uncertainty by revising decisions throughout a multistage planning horizon based on the latest information. 
Although this approach has been widely applied in different application domains (e.g., \citet{gul2015progressive, pan2016strong}), 
its widespread adoption in practice can be challenging for various problems 
due to the limited flexibility required within the decision-making processes, preventing revising decisions at each stage. To address this challenge, this paper advances multistage stochastic programming by incorporating the need for limited flexibility and commitments through optimally determining the most critical stages \textit{when} decisions can be revised.  

We introduce a new optimization approach called  \textit{adaptive multistage stochastic programming}, which limits the number of revisions allowed throughout the planning horizon for a subset of decisions that require a certain level of commitment. In the first stage, we determine the optimal stages of revisions for these decisions. Subsequently, at these selected stages, these decisions can be revised similarly to classic multistage stochastic programming, while in the remaining stages, commitment is ensured. The clear advantage of our approach is that at the start of the planning horizon, it is already known when decisions will be revised, which allows for effective communication between dependent actors in decision-making contexts.
In the remainder of this section, we discuss related work and summarize our contributions by introducing adaptive multistage stochastic programming.

\subsection{Related Work} \label{sec:lit}

Stochastic programming is one of the most effective analytical techniques for decision-making under uncertainty \citep{birge2000introduction}. It involves modeling the stochastic process using scenario trees \citep{shapiro2009lectures}, which represent the possible uncertainty realizations at different stages. There are two fundamental techniques in stochastic programming: two-stage and multistage, which differ in how they handle sequential decision-making under uncertainty. These techniques generally involve two types of decisions: ``here-and-now" decisions that can be immediately executed, and ``wait-and-see" decisions that depend on the system's state (i.e., the observation of the system in each decision stage), influenced by prior decisions and observations of uncertain parameters. While two-stage stochastic programming provides fixed decisions upfront planning for here-and-now decisions, multistage stochastic programming extends the framework by revising all decisions at each stage based on the uncertainty realized so far \citep{birge1985decomposition}. 


One critical type of uncertainty that arises in stochastic optimization is endogenous uncertainty, which originates from the decisions of the decision-maker. Endogenous uncertainty adds an additional layer of complexity to the associated optimization problems, as the decision-maker must consider not only the uncertain conditions they face but also the potential impact of their own decisions on the underlying stochastic process in future decision stages.
Most often, the stochastic programming literature considers exogenous uncertainty, where decisions do not impact the underlying stochastic process. Relatively limited attention has been paid to endogenous uncertainty. In this type of uncertainty, decisions can change the uncertainty realizations by affecting the underlying probability distributions \citep[see, e.g., ][]{hellemo2018decision,basciftci2021distributionally,csafak2022two,drusvyatskiy2022stochastic}
or alternatively the time which uncertainty realization is observed \citep[see, e.g.,][]{goel2006class, gupta2011solution, apap2017models}, where the latter setting will be discussed further. 

The endogenous uncertainty literature involving this second line of research allows revising decisions in each stage based on the observed uncertainty. 
However, it is important to note that uncertainty is not automatically realized in each stage, and its realization depends on the specific decisions made. Until these decisions are made, it remains unknown which uncertain parameters will be observed, and there is a possibility that some uncertain parameters may not be observed at all.
To address this, non-anticipativity constraints (NACs) are commonly used for restricting the feasible set of decision variables depending on the time of observed uncertainty, where the number of NACs exponentially increases with the number of possible uncertainty realizations. \cite{goel2006class} present a mixed-integer disjunctive formulation for stochastic programming that enforces endogenous uncertainty with conditional NACs and investigates theoretical properties to eliminate a particular set of NACs. Despite the elimination techniques, solving the reduced models becomes challenging as the model size increases.  For example, given that a small number of NACs are active at the optimal solution, \cite{colvin2010modeling} develop a branch and cut algorithm in which a certain set of NACs are initially removed and reinserted if they are violated within the search tree. In a subsequent study,  \cite{gupta2011solution} propose NACs elimination together with heuristic solution strategies: including NACs only up to a limited number of stages, 
and using Lagrangean relaxation of the NACs to decompose the problem into scenarios, which is later extended to the settings with both endogenous and exogenous uncertainties \citep{apap2017models}. 
To reduce the number of NACs, \cite{boland2016minimum} provide general conditions depending on the characteristics of the scenario structure. 
As another relevant line of research, \cite{gupta2016algorithms} and \cite{singla2018combinatorial} study stochastic probing algorithms for specific problems, determining uncertainty realization through probing decisions and associated costs. Although those studies show similarities as the decisions affect future conditions and their setting requires NACs, they differ from our study by making the uncertainty realization conditional and consequently not fully realizing the underlying stochastic process. 
Furthermore, existing studies do not focus on providing a generic optimization under uncertainty framework under partial flexibility, and present different formulation techniques and NACs reductions that are not directly applicable to our  proposed setting. 

This paper introduces Adaptive Multistage Stochastic Programming; a new optimization approach that differs from the existing literature by leveraging decision commitment into multistage stochastic programming. It determines at the start of the planning horizon \textit{when} revisions can be made subject to a limited number of allowed revisions. \cite{basciftci2019adaptive} propose an approach called adaptive two-stage stochastic programming in which they optimize the timing of a single decision revision. They propose analytical approaches over a capacity expansion planning problem through a bound analysis of the adaptive two-stage stochastic program in comparison to two-stage and multistage stochastic programming alternatives, based on the timing of the revision decision. Although this is the first paper to consider revisions as decision variables, we observe two restrictions in their model: i) their model allows revising the decision only once, which is restrictive if the planning horizon becomes larger, and ii) they do not provide elimination techniques for the NACs that would reduce the computational burden of the problem. Besides, their study does not provide a generic solution algorithm, but rather approximation algorithms tailored to address the capacity expansion planning problem. 

Even without the NACs needed for Adaptive Multistage Stochastic Programming, solving a classic multistage stochastic programming problem is theoretically and computationally challenging. Furthermore, as the number of stages increases, the size of the scenario tree grows exponentially, making it even more difficult to solve \citep{shapiro2005complexity}. Decomposition methods are commonly used to solve large-scale optimization models, and two well-known approaches are the integer L-shaped method \citep{laporte1993integer,angulo2016improving} and Benders decomposition  \citep{angulo2016improving,fischetti2017redesigning,rahmaniani2020benders}, which can be further utilized over stochastic programs. The L-shaped method is a cutting plane technique that decomposes a model into master problem and subproblems to separate complicating variables from the remainder of the model \citep{van1969shaped}. It was originally developed for continuous subproblems and later improved to the integer L-shaped method \citep{laporte1993integer}, which handles pure binary decisions in the master problem and mixed-integer decisions in the subproblem. Benders decomposition is another approach that solves models with mixed-integer master problems and linear subproblems \citep{benders}. We provide a solution approach for generic adaptive multistage stochastic programs with mixed-integer variables by adapting these decomposition methods and leveraging the structure of our proposed formulation. 

\subsection{Contributions} 
Summarizing, we make the following scientific contributions.
\begin{itemize}
    \item We introduce adaptive multistage stochastic programming: a novel and generic decision-making approach that incorporates a balance between decision flexibility and commitment into multistage stochastic programming. By considering the flexibility level of the decision-making process, we set a maximum limit on the number of revisions. The main factor in our optimization problem is determining the optimal stages for revising each decision. While these optimal stages are determined in the first stage, the corresponding revisions are made later at those optimal stages. Furthermore, we prove the NP-hardness of adaptive multistage stochastic programming and establish its connection to classic models such as two-stage stochastic programming and multistage stochastic programming. 
    \item To address the complexity of NACs in adaptive multistage stochastic programming problems, we present theoretical properties that reduce the problem size significantly. In addition, we propose a tailored cutting plane-based solution algorithm that customizes the integer L-shaped method and Benders decomposition to be compatible with mixed-integer master problem and subproblems. Furthermore, we provide a preprocessing approach that eliminates a set of candidate feasible solutions through a cut-generation procedure. Our findings show that these strategies significantly enhance the computational efficiency of decision-making processes.
    \item We conduct computational experiments on adaptive stochastic programming variants of two classical optimization problems; stochastic lot-sizing and generation expansion planning. 
    Despite the limited flexibility, we show that it is possible to converge to the objective of multistage stochastic programming if revisions are made at optimal stages. It is remarkable that optimal stage selection in less flexible case can be more efficient than arbitrary selection in more flexible case. 
    We conclude that adaptive multistage stochastic programming offers businesses practically feasible solutions that ensure decision commitment in their decision-making process while staying flexible to act upon realized uncertainty at the most critical stages.
    \item The generation expansion problem is calibrated with real-world data. Through this analysis, we provide managerial insights on optimal revision stages depending on cost and technical parameters of generation sources. 
    We further show the computational efficiency of our reformulations and solution algorithms, where our NACs elimination techniques reduce computation times by 59\%, and our solution algorithm further reduces it by another 56\% on average considering various instances.
\end{itemize}

The remainder of the paper is organized as follows. Section \ref{sec:Model} introduces our Adaptive Multistage Stochastic Programming approach. Section \ref{sec:NACsElimination} presents a series of NACs elimination techniques to enhance the proposed model. Section \ref{sec:SolutionApproach} provides the solution approach. Afterward, Section \ref{sec:ComputationalExperiments} shows the value of adaptive multistage stochastic programming by computational experiments on two canonical optimization problems; stochastic lot-sizing and generation expansion planning with managerial insights. 
Section \ref{sec:Conclusions} concludes our paper with final remarks. 

\section{Adaptive Multistage Stochastic Programming }\label{sec:Model}
In this section, we introduce Adaptive Multistage Stochastic Programming (AMSP) for sequential decision-making problems under uncertainty where an optimal balance between decision flexibility and commitment is required. AMSP limits the number of revisions throughout the planning horizon and optimally determines stages at which decision revisions will be allowed for a subset of decisions requiring partial flexibility. Consequently, these decisions can be revised only at those selected stages. 

To clarify this new approach, we first present multistage and two-stage stochastic programming formulations, which allow for different flexibility levels for these decisions. In particular, multistage stochastic programming offers a fully flexible setting that allows for revisions as new information is revealed through the planning horizon. In contrast, the two-stage stochastic programming commits a subset of decisions, i.e., here-and-now decisions, statically upfront in planning and does not allow revising these decisions. To this end, AMSP can be interpreted as a more generic approach since it captures both multistage and two-stage stochastic programming depending on the flexibility level allowed in decision-making. 
Combining the above, we present our AMSP approach in the following order. We first define the underlying stochastic process based on a scenario tree, which is a fundamental approach for representing uncertainty in sequential decision-making problems \citep{Ruszczynski2003}. We then provide two generic formulations for multistage and two-stage stochastic programming and show the associated decision structures on how decisions are revised over the scenario tree. Next, we introduce our AMSP formulation and show its decision structure in comparison with two-stage and multistage stochastic programming. Lastly, we discuss several theoretical properties and prove NP-hardness of AMSP.

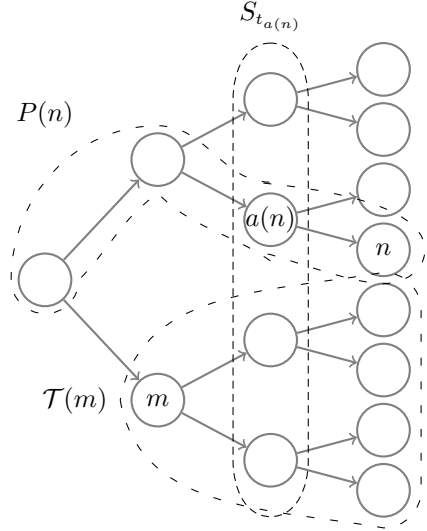
\begin{figure}[h]
    \centering
        \caption{Scenario Tree with $T=4$ Stages and $B=2$ Branches}
    \label{fig:Scenario Tree}
\begin{tikzpicture}
{\color{gray}
 \tikzstyle{annot} = [text width=4em, text centered]
\begin{scope}[every node/.style={circle,thick,draw,minimum size=0.7cm}]
\node (I0) at (0,0) {};
\node (I1) at (1.5,1.6) {};
\node (I2) at (1.5,-1.6) {{\color{black}\small{$m$}}};
\node (I3) at (3,2.4) {};
\node (I4) at (3,0.8) {};
\node (I5) at (3,-0.8) {};
\node (I6) at (3,-2.4) {};
\node (I7) at (4.5,2.8) {};
\node (I8) at (4.5,2) {};
\node (I9) at (4.5,1.2) {};
\node (I10) at (4.5,0.4) {{\color{black}\small{$n$}}};
\node (I11) at (4.5,-0.4) {};
\node (I12) at (4.5,-1.2) {};
\node (I13) at (4.5,-2) {};
\node (I14) at (4.5,-2.8) {};
\end{scope}
}
\begin{scope}[every node/.style={fill=white,circle},
              every edge/.style={draw=gray,thick}]
     \path [->] (I0) edge (I1);
            \path [->] (I0) edge (I2);
            \path [->] (I1) edge (I3);
            \path [->] (I1) edge (I4);
            \path [->] (I2) edge (I5);
            \path [->] (I2) edge (I6);
            \path [->] (I3) edge (I7);
            \path [->] (I3) edge (I8);
            \path [->] (I4) edge (I9);
            \path [->] (I4) edge (I10);
            \path [->] (I5) edge (I11);
            \path [->] (I5) edge (I12);
            \path [->] (I6) edge (I13);
            \path [->] (I6) edge (I14);
\end{scope}

 \begin{scope}[fill opacity=0, densely dashed]
    \draw ($(I3)+(0,0.75)$)
        to[out=180,in=90] ($(I3) + (-0.5,0)$)
        to[out=270,in=90] ($(I6) + (-0.5,0)$)  
        to[out=270,in=180] ($(I6) + (0,-0.75)$)
        to[out=0,in=270] ($(I6) + (0.5,0)$)
        to[out=90,in=270] ($(I3) + (0.5,0)$)  
        to[out=90,in=0] ($(I3) + (0,0.75)$)   ;
        \end{scope}

\begin{scope}[fill opacity=0, loosely dashed]
    \draw ($(I2)+(-0.5,0)$)
        to[out=280,in=170] ($(I14) + (0,-0.5)$)  
        to[out=0,in=270] ($(I14) + (0.5,0)$)
        to[out=90,in=270] ($(I11) + (0.5,0)$)
        to[out=90,in=0] ($(I11) + (0,0.5)$)  
        to[out=190,in=80] ($(I2) + (-0.5,0)$)   ;
        \end{scope}

\begin{scope}[fill opacity=0, loosely dashed]
 \draw ($(I0)+(-0.45,0)$)
    to[out=90,in=180] ($(I1)+(0,0.45)$)
    to[out=0,in=180] ($(I4)+(0,0.5)$)
    to[out=180,in=90] ($(I10)+(0.55,0)$)
    to[out=270,in=0] ($(I10)+(0,-0.45)$)
    to[out=180,in=145] ($(I4)+(0,-0.55)$)
    to[out=145,in=315] ($(I1)+(0,-0.45)$)
    to[out=215,in=45] ($(I0)+(0.45,0.0)$)
    to[out=270,in=0] ($(I0)+(0.0,-0.45)$)
    to[out=180,in=270] ($(I0)+(-0.45,0.0)$) ;
        \end{scope}
       
\node[above of=I3, node distance=1.1cm] (hl) {\small{$S_{t_{a(n)}}$}};
\node[above of=I0, node distance=2.2cm] (h) {\small{$P(n)$}};  
\node[above of=I4, node distance=0cm] (h) {\small{$a(n)$}}; 
\node[left of=I2, node distance=1.1cm] (h) {\small{$\mathcal{T}(m)$}}; 
\end{tikzpicture} 
\end{figure}

We consider a scenario tree $\mathcal{T}$ to model the uncertainty structure throughout planning horizon of $T$ stages. We represent each node of the scenario tree as $n \in \mathcal{T}$, and each stage as $t \in \mathcal{P} := \{1, \dots, T\}$. We define the set of nodes in each stage $t \in \mathcal{P}$ as $S_{t}$, and the stage of a node $n$ as $t_{n} \in \mathcal{P}$. Each node in the scenario tree, except for the leaf nodes, has $B$ branches that represent possible outcomes in the subsequent stage. 
Each node $n$, except the root node, has an ancestor denoted by $a(n)$. The unique path from the root node to a specific node $n$ is represented by $P(n)$. Each path from the root node to a leaf node corresponds to a scenario, indicating that every $P(n)$ represents a scenario when $n \in S_{T}$. We denote the sub-tree rooted at node $n$ until stage $t$ as $\mathcal{T}(n, t)$ for $t_{n} \leq t \leq T$. To shorten the notation, when the last stage of the sub-tree is $T$, we let $\mathcal{T}(n):=\mathcal{T}(n, T)$ for all $n \in \mathcal{T}$. The probability of observing each node $n$ is given by $p_{n}$, where $\sum_{n \in S_{t}} p_{n}=1$ for all $t \in \mathcal{P}$. Figure~\ref{fig:Scenario Tree} shows an example scenario tree to illustrate these concepts.

The decision variables corresponding to node $n \in \mathcal{T}$ consist of the \textit{state} variables $\left\{\mathbf{x_{n}}\right\}_{n \in \mathcal{T}}$, and the \textit{stage} variables $\left\{\mathbf{y_{n}}\right\}_{n \in \mathcal{T}}$. In particular, state variables carry information between different stages, whereas the stage variables at node $n \in \mathcal{T}$ are local variables to their associated stage $t_n \in \mathcal{P}$. We let the dimensions of the variables $\mathbf{x_n} = (x_{1n},\dots,x_{In})$ and  $\mathbf{y_n} = (y_{1n},\dots,y_{Jn})$  be $I$ and $J$ for all $n \in \mathcal{T}$, respectively, and define the sets $\mathcal{I} = \{1,\dots,I\}$, and $\mathcal{J} = \{1,\dots,J\}$. 

Here, we consider a generic formulation where these variables can be integer or continuous. The parameters of each node $n \in \mathcal{T}$ are represented as $\left(\mathbf{a_{n}}, \mathbf{b_{n}}, \mathbf{C_{nm}}, \mathbf{D_{n}}, \mathbf{d_{n}}\right)$. For a given scenario tree, we define multistage stochastic programming (MSP), where both state and stage decisions can be revised in each stage based on the observed uncertainty, as follows:
\begin{subequations}\label{MSLP}
\allowdisplaybreaks\begin{align}
(MSP) \ \min_{\mathbf{x}, \mathbf{y}} \quad & \sum_{n \in \mathcal{T}} p_{n}\left(\mathbf{a_{n}}^{\top} \mathbf{x_{n}}+ \mathbf{b_{n}}^{\top} \mathbf{y_{n}}\right) \label{eq:MSObj}\\
\text{s.t.} \quad &  \sum_{m \in P(n)} \mathbf{C_{n m} x_{m}}+ \mathbf{D_{n} y_{n}} \geq \mathbf{d_{n}} \quad \forall n \in \mathcal{T} \label{eq:MScons}\\
&  \mathbf{x_{n}}  \in \mathcal{X}_{n}, \mathbf{y_{n}} \in \mathcal{Y}_{n}, \quad \forall n \in \mathcal{T} \label{eq:MSdom}
\end{align}
\end{subequations}
Here, objective \eqref{eq:MSObj} scales the cost $\mathbf{a_n}$ and $\mathbf{b_n}$ of $\mathbf{x_n}$ and $\mathbf{y_n}$, respectively, with the probability $p_n$ of observing node $n \in \mathcal{T}$. Constraints \eqref{eq:MScons} link state and stage variables, and constraints \eqref{eq:MSdom} impose local constraints over each node $n \in \mathcal{T}$ on the variables $\mathbf{x_n}$ and $\mathbf{y_n}$ by the sets $\mathcal{X}_{n}$ and $\mathcal{Y}_{n}$. In formulation $(MSP)$, we allow state variables $\mathbf{x_n}$ to differ in each node $n \in \cal T$. 

An often-used approximation to multistage stochastic programming, called two-stage stochastic programming, results from imposing equality of the state variables $\mathbf{x_n}$ in each stage $t_n \in \mathcal{P}$. This restriction implies that state decisions cannot adapt to realized uncertainty,  making it less flexible than formulation $(MSP)$. Note that,  in two-stage stochastic programming, uncertain parameters are typically defined over scenarios which are represented as paths within the scenario tree. However, as in \cite{huang2009value}, we define these parameters at the node level instead of the scenario level. The formulation of the two-stage stochastic programming is given as follows:
\begin{subequations}\label{TSLP}
\allowdisplaybreaks\begin{align}
(2SP) \ \min_{\mathbf{x}, \mathbf{y}} \quad & \sum_{n \in \mathcal{T}} p_{n}\left(\mathbf{a_{n}}^{\top} \mathbf{x_{n}}+\mathbf{b_{n}}^{\top} \mathbf{y_{n}}\right) \label{eq:2sObj}\\
\text{s.t.} \quad & \eqref{eq:MScons}, \eqref{eq:MSdom} \notag \\
& \mathbf{x_m}=\mathbf{x_n} \quad  \forall m, n \in S_{t}, t \in \mathcal{P}\label{Two-stageC}
\end{align}
\end{subequations}
In formulation $(2SP)$, objective \eqref{eq:2sObj}, and constraints \eqref{eq:MScons} and \eqref{eq:MSdom} are the same as in formulation $(MSP)$, whereas, with constraints \eqref{Two-stageC}, we impose equality of state variables across different nodes of the same stage, which can be referred as NACs. Before introducing AMSP, we demonstrate the differences in decision structures between formulations $(MSP)$ and $(2SP)$.

\begin{figure}[h]
\centering
\caption{Decision structures of state and stage variables}
\vspace{2mm}
\label{fig:Decisions_multi_two}
\begin{subfigure}{.26\textwidth}
\captionsetup{size=small}
  \caption{ \textit{State} variables in (MSP)}
  \label{fig:multi_stage}
  \begin{tikzpicture}[scale=0.8]
 \tikzstyle{annot} = [text width=4em, text centered]
\begin{scope}[every node/.style={circle,thick,draw,minimum size=0.7cm,scale=0.7}]
\node  (I0) at (0,0) {};
\node  (I1) at (1.25,1.6) {};
\node (I2) at (1.25,-1.6) {};
\node  (I3) at (2.5,2.4) {};
\node  (I4) at (2.5,0.8) {};
\node  (I5) at (2.5,-0.8) {};
\node  (I6) at (2.5,-2.4) {};
\node (I7) at (3.75,2.8) {};
\node  (I8) at (3.75,2) {};
\node  (I9) at (3.75,1.2) {};
\node  (I10) at (3.75,0.4) {};
\node (I11) at (3.75,-0.4) {};
\node (I12) at (3.75,-1.2) {};
\node (I13) at (3.75,-2) {};
\node (I14) at (3.75,-2.8) {};
\end{scope}
\begin{scope}[every node/.style={fill=white,circle},
              every edge/.style={draw=black,thick}]
     \path [->] (I0) edge (I1);
            \path [->] (I0) edge (I2);
            \path [->] (I1) edge (I3);
            \path [->] (I1) edge (I4);
            \path [->] (I2) edge (I5);
            \path [->] (I2) edge (I6);
            \path [->] (I3) edge (I7);
            \path [->] (I3) edge (I8);
            \path [->] (I4) edge (I9);
            \path [->] (I4) edge (I10);
            \path [->] (I5) edge (I11);
            \path [->] (I5) edge (I12);
            \path [->] (I6) edge (I13);
            \path [->] (I6) edge (I14);
\end{scope}
\draw [dashed,white] (0,3.1) -- (0,3.05);
\draw [dashed,white] (0,-3.1) -- (0,-3.05);
\end{tikzpicture} 
\end{subfigure}%
\hspace{0.2cm}
\begin{subfigure}{.25\textwidth}
\captionsetup{size=small}
  \caption{\textit{State} variables in (2SP)}
  \label{fig:two_stage}
  \begin{tikzpicture}[scale=0.8]
 \tikzstyle{annot} = [text width=4em, text centered]
\begin{scope}[every node/.style={circle,thick,draw,minimum size=0.7cm,scale=0.7}]
\node  (I0) at (0,0) {};
\node  (I1) at (1.25,0) {};
\node  (I3) at (2.5,0) {};
\node (I7) at (3.75,0) {};
\end{scope}
\begin{scope}[every node/.style={fill=white,circle},
              every edge/.style={draw=black,thick}]
     \path [->] (I0) edge (I1);
            \path [->] (I1) edge (I3);
            \path [->] (I3) edge (I7);
\end{scope}
\draw [dashed,white] (0,3.1) -- (0,3.05);
\draw [dashed,white] (0,-3.1) -- (0,-3.05);
\end{tikzpicture} 
\end{subfigure}%
\hspace{0.2cm}
\begin{subfigure}{.35\textwidth}
\captionsetup{size=small}
  \caption{\textit{Stage} variables in (2SP) \& (MSP)}
  \label{fig:stage variables}
  \begin{tikzpicture}[scale=0.8]
 \tikzstyle{annot} = [text width=4em, text centered]
\begin{scope}[every node/.style={circle,thick,draw,minimum size=0.7cm,scale=0.7}]
\node  (I0) at (0,0) {};
\node  (I1) at (1.25,1.6) {};
\node (I2) at (1.25,-1.6) {};
\node  (I3) at (2.5,2.4) {};
\node  (I4) at (2.5,0.8) {};
\node  (I5) at (2.5,-0.8) {};
\node  (I6) at (2.5,-2.4) {};
\node (I7) at (3.75,2.8) {};
\node  (I8) at (3.75,2) {};
\node  (I9) at (3.75,1.2) {};
\node  (I10) at (3.75,0.4) {};
\node (I11) at (3.75,-0.4) {};
\node (I12) at (3.75,-1.2) {};
\node (I13) at (3.75,-2) {};
\node (I14) at (3.75,-2.8) {};
\end{scope}

\begin{scope}[every node/.style={fill=white,circle},
              every edge/.style={draw=black,thick}]
     \path [->] (I0) edge (I1);
            \path [->] (I0) edge (I2);
            \path [->] (I1) edge (I3);
            \path [->] (I1) edge (I4);
            \path [->] (I2) edge (I5);
            \path [->] (I2) edge (I6);
            \path [->] (I3) edge (I7);
            \path [->] (I3) edge (I8);
            \path [->] (I4) edge (I9);
            \path [->] (I4) edge (I10);
            \path [->] (I5) edge (I11);
            \path [->] (I5) edge (I12);
            \path [->] (I6) edge (I13);
            \path [->] (I6) edge (I14);
\end{scope}
\draw [dashed,white] (0,3.1) -- (0,3.05);
\draw [dashed,white] (0,-3.1) -- (0,-3.05);
\end{tikzpicture} 
\end{subfigure}
\end{figure}
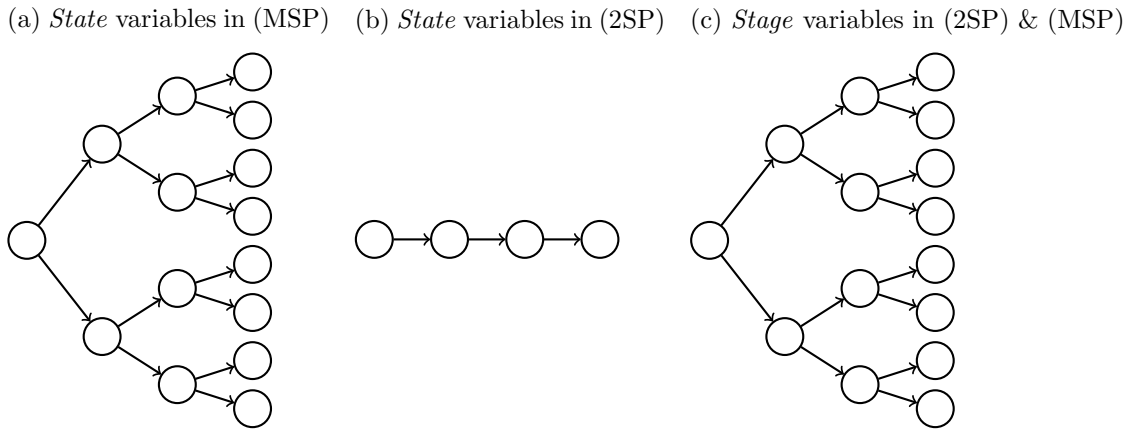
\vspace{-8mm}
\begin{example}[Decision structures in formulations $(MSP)$ and $(2SP)$]
We illustrate decision structures over 
state variables $\{\mathbf{x_n}\}_{n \in \mathcal{T}}$
and stage variables $\{\mathbf{y_n}\}_{n \in \mathcal{T}}$ over a scenario tree with four stages and two branches (e.g., the scenario tree in Figure~\ref{fig:Scenario Tree}) in Figure~\ref{fig:Decisions_multi_two}. In multistage stochastic programming, in Figure~\ref{fig:multi_stage}, state variables can be determined for each node $n \in \mathcal{T}$ as uncertainty is realized over the planning horizon. In two-stage stochastic programming, in Figure~\ref{fig:two_stage}, the values of the state variables for each stage are determined at the beginning of the planning horizon for $t = 1, \dots, T$, and cannot be adjusted for each realization. However, the decision structure of the stage variables $\{\mathbf{y_n}\}_{n \in \mathcal{T}}$ are the same under both approaches, and they can be determined as uncertainty is revealed in each stage as depicted in Figure~\ref{fig:stage variables}.  
\end{example}

A decision-maker can desire a balance between decision flexibility and commitment, which requires more flexibility than formulation $(2SP)$ and less flexibility than formulation $(MSP)$ in revising state variables. To address this need, AMSP limits the change of each state variable over the multi-period planning horizon and includes the decisions on \textit{when} the state variables are allowed to be updated within planning. We refer to the stage at which a state variable can be changed based on uncertainty observation as \textit{revision point}, and the changes themselves as \textit{revisions}. Crucial for AMSP is that the revision points are part of the optimization problem as here-and-now decisions and decided before uncertainties are realized. State variables can only be revised at these revision points, while stage variables can be revised at each stage. Another notable advantage of AMSP is its capability to provide revision points specific to each state variable.



To ensure commitment over state decisions, we limit flexibility by allowing at most $\mu$ revision points, where $\mu \in \{0, 1,\dots, T-1\}$. 
We define a set of revision points for each state variable, which are denoted by $t^*_{ik} \in \{2, \dots, T \}$, representing the $k^{th}$ revision time of variable $\mathbf{x_i}$ for every $i \in \mathcal{I}$ and $k \in \{1,\dots,\mu\}$,  
which are optimally determined by the adaptive multistage program. We note that $t^*_{ik} < t^*_{i(k+1)}$ should be satisfied for every $k = {1,\dots,\mu-1}$, as $(k+1)^{st}$ revision can only occur after $k^{th}$ revision takes place. 
Before introducing the formulation of AMSP, we illustrate various decision structures of the state variables under different revision points. 


\begin{figure}[h]
\centering
\caption{Decision structures of state variable $i$ for AMSP}
\label{fig:Decisions_adaptive}
\begin{subfigure}{.25\textwidth}
\captionsetup{size=small}
  \caption{$\mu=1, t_{i1}^*=2$}
  \label{fig:mu_1_t_2}
  \begin{tikzpicture}[scale=0.8]
 \tikzstyle{annot} = [text width=4em, text centered]
\begin{scope}[every node/.style={circle,thick,draw,minimum size=0.7cm,scale=0.7}]
\node (I0) at (0,0) {};
\node (I1) at (1.25,1.6) {};
\node (I2) at (1.25,-1.6) {};
\node (I3) at (2.5,1.6) {};
\node (I5) at (2.5,-1.6) {};
\node (I7) at (3.75,1.6) {};
\node (I11) at (3.75,-1.6) {};
\end{scope}
\begin{scope}[every node/.style={fill=white,circle},
              every edge/.style={draw=black,thick}]
     \path [->] (I0) edge (I1);
            \path [->] (I0) edge (I2);
            \path [->] (I1) edge (I3);
            \path [->] (I2) edge (I5);
            \path [->] (I3) edge (I7);
            \path [->] (I5) edge (I11);

\end{scope}
\draw [dashed] (1.25,2.8) -- (1.25,-2.8);
\draw [dashed,white] (0,3.1) -- (0,3.05);
\draw [dashed,white] (0,-3.1) -- (0,-3.05);
\end{tikzpicture} 
\end{subfigure}%
\begin{subfigure}{.25\textwidth}
\captionsetup{size=small}
\caption{$\mu=1, t_{i1}^*=4$}
  \label{fig:mu_1_t_4}
  \begin{tikzpicture}[scale=0.8]
 \tikzstyle{annot} = [text width=4em, text centered]
\begin{scope}[every node/.style={circle,thick,draw,minimum size=0.7cm,scale=0.7}]
\node  (I0) at (0,0) {};
\node  (I1) at (1.25,0) {};
\node  (I3) at (2.5,0) {};
\node (I7) at (3.75,2.8) {};
\node  (I8) at (3.75,2) {};
\node  (I9) at (3.75,1.2) {};
\node  (I10) at (3.75,0.4) {};
\node (I11) at (3.75,-0.4) {};
\node (I12) at (3.75,-1.2) {};
\node (I13) at (3.75,-2) {};
\node (I14) at (3.75,-2.8) {};
\end{scope}

\begin{scope}[every node/.style={fill=white,circle},
              every edge/.style={draw=black,thick}]
     \path [->] (I0) edge (I1);
            \path [->] (I1) edge (I3);
            \path [->] (I3) edge (I7);
            \path [->] (I3) edge (I8);
            \path [->] (I3) edge (I9);
            \path [->] (I3) edge (I10);
            \path [->] (I3) edge (I11);
            \path [->] (I3) edge (I12);
            \path [->] (I3) edge (I13);
            \path [->] (I3) edge (I14);
\end{scope}
\draw [dashed] (3.75,2.8) -- (3.75,-2.8);
\draw [dashed,white] (0,3.1) -- (0,3.05);
\draw [dashed,white] (0,-3.1) -- (0,-3.05);
\end{tikzpicture} 
\end{subfigure}%
\begin{subfigure}{.25\textwidth}
\captionsetup{size=small}
  \caption{$\mu=2, t_{i1}^*=2, t_{i2}^*=3$}
  \label{fig:mu_2_t_2_3}
  \begin{tikzpicture}[scale=0.8]
 \tikzstyle{annot} = [text width=4em, text centered]
\begin{scope}[every node/.style={circle,thick,draw,minimum size=0.7cm,scale=0.7}]
\node  (I0) at (0,0) {};
\node  (I1) at (1.25,1.6) {};
\node (I2) at (1.25,-1.6) {};
\node  (I3) at (2.5,2.4) {};
\node  (I4) at (2.5,0.8) {};
\node  (I5) at (2.5,-0.8) {};
\node  (I6) at (2.5,-2.4) {};
\node (I7) at (3.75,2.4) {};
\node  (I9) at (3.75,0.8) {};
\node (I11) at (3.75,-0.8) {};
\node (I13) at (3.75,-2.4) {};
\end{scope}

\begin{scope}[every node/.style={fill=white,circle},
              every edge/.style={draw=black,thick}]
     \path [->] (I0) edge (I1);
            \path [->] (I0) edge (I2);
            \path [->] (I1) edge (I3);
            \path [->] (I1) edge (I4);
            \path [->] (I2) edge (I5);
            \path [->] (I2) edge (I6);
            \path [->] (I3) edge (I7);
            \path [->] (I4) edge (I9);     
            \path [->] (I5) edge (I11);
            \path [->] (I6) edge (I13);

\end{scope}
\draw [dashed] (1.25,2.8) -- (1.25,-2.8);
\draw [dashed] (2.5,2.8) -- (2.5,-2.8);
\draw [dashed,white] (0,3.1) -- (0,3.05);
\draw [dashed,white] (0,-3.1) -- (0,-3.05);
\end{tikzpicture} 
\end{subfigure}%
\begin{subfigure}{.25\textwidth}
\captionsetup{size=small}
  \caption{$\mu=2, t_{i1}^*=2, t_{i2}^*=4$}
  \label{fig:mu_2_t_2_4}
 \begin{tikzpicture}[scale=0.8]
 \tikzstyle{annot} = [text width=4em, text centered]
\begin{scope}[every node/.style={circle,thick,draw,minimum size=0.7cm,scale=0.7}]
\node  (I0) at (0,0) {};
\node  (I1) at (1.25,1.6) {};
\node (I2) at (1.25,-1.6) {};
\node  (I3) at (2.5,1.6) {};
\node  (I5) at (2.5,-1.6) {};
\node (I7) at (3.75,2.8) {};
\node  (I8) at (3.75,2) {};
\node  (I9) at (3.75,1.2) {};
\node  (I10) at (3.75,0.4) {};
\node (I11) at (3.75,-0.4) {};
\node (I12) at (3.75,-1.2) {};
\node (I13) at (3.75,-2) {};
\node (I14) at (3.75,-2.8) {};
\end{scope}
\begin{scope}[every node/.style={fill=white,circle},
              every edge/.style={draw=black,thick}]
     \path [->] (I0) edge (I1);
            \path [->] (I0) edge (I2);
            \path [->] (I1) edge (I3);
            \path [->] (I2) edge (I5);
            \path [->] (I3) edge (I7);
            \path [->] (I3) edge (I8);
            \path [->] (I3) edge (I9);
            \path [->] (I3) edge (I10);
            \path [->] (I5) edge (I11);
            \path [->] (I5) edge (I12);
            \path [->] (I5) edge (I13);
            \path [->] (I5) edge (I14);
\end{scope}
\draw [dashed] (1.25,2.8) -- (1.25,-2.8);
\draw [dashed] (3.75,2.8) -- (3.75,-2.8);
\draw [dashed,white] (0,3.1) -- (0,3.05);
\draw [dashed,white] (0,-3.1) -- (0,-3.05);
\end{tikzpicture} 
\end{subfigure}
\end{figure}
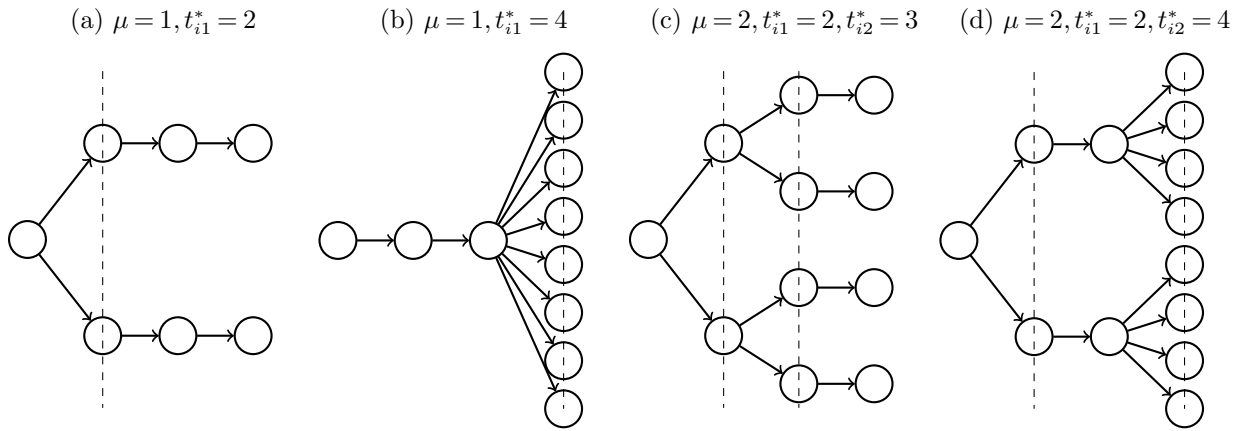

\vspace{-8mm}
\begin{example}[Decision structures of state variables in AMSP]
We highlight how AMSP works by showing four examples of decision structures under given revision point decisions in Figure~\ref{fig:Decisions_adaptive}. We consider a state variable $\{x_{in}\}_{n \in \mathcal{T}}$, for a given $i$, considering the scenario tree in Figure~\ref{fig:Scenario Tree}. In Figure~\ref{fig:mu_1_t_2} and \ref{fig:mu_1_t_4}, $\mu$ is set to $1$, which implies that only one revision is allowed for the state variable  $\{x_{in}\}_{n \in \mathcal{T}}$. In Figure~\ref{fig:mu_1_t_2}, the second stage is determined as the revision stage, and decision variables corresponding to the second stage's nodes can differ based on the realized uncertainty in stage two. However, the nodes for the third and fourth stages at each sub-tree rooting from nodes of the second stage are condensed, meaning that a single decision is made for those condensed nodes.  
In other words, one can commit to decisions made for third and fourth stages in the second stage. In Figure~\ref{fig:mu_1_t_4}, the revision stage is the fourth stage. Thus, nodes at the first three stages are condensed until this time. Furthermore, Figures~\ref{fig:mu_2_t_2_3} and \ref{fig:mu_2_t_2_4} present cases where $\mu$ is set to 2, where the decisions are revised in the second and third stages and the second and fourth stages, respectively. 
\end{example}

We formulate AMSP as a stochastic program considering $\mu$ revision points for each state variable over the planning horizon. We include revision points $\Gamma^*_i = \{t^*_{i1},\dots,t^*_{i\mu}\}, i \in \mathcal{I}$ as decision variables 
by imposing equalities between the state variables depending on these revision points, leading to the following formulation: 
\begin{subequations}\label{AMSLP}
\begin{align}
 \min_{\mathbf{x}, \mathbf{y}, \mathbf{t^*}} \quad & \sum_{n \in \mathcal{T}} p_{n}\left(\mathbf{a_{n}}^{\top} \mathbf{x_{n}}+  \mathbf{b_{n}}^{\top} \mathbf{y_{n}}\right) \label{AMSLP_obj} \\
\text{s.t.} \quad &  \eqref{eq:MScons}, \eqref{eq:MSdom} \notag \\ 
& x_{im}=x_{in} \quad \forall m, n \in S_t, \quad t<t^*_{i1}, \quad i \in \mathcal{I} \label{AMSLP_d}\\
& x_{im}=x_{in} \quad \forall m, n \in S_t \cap \mathcal{T}(l,t^*_{i(k+1)}), \quad l \in S_{ t^*_{ik}}, \quad  t^*_{i(k+1)} >  t \geq t^*_{ik}, \quad \mu-1 \geq k \geq 1, \quad i \in \mathcal{I}\label{AMSLP_e} \\
& x_{im}=x_{in} \quad \forall m, n \in S_t \cap \mathcal{T}(l), \quad l \in S_{t^*_{i\mu}}, \quad  t \geq t^*_{i\mu}, \quad i \in \mathcal{I}\label{AMSLP_f} \\ 
& t^*_{ik} \in \Gamma_i^*  \subseteq  \{2, \dots, T\}  \quad \forall i \in \mathcal{I}, \quad k \in \{1,\dots,\mu\} \label{AMSLP_g} 
\end{align}
\end{subequations}
Objective \eqref{AMSLP_obj}, constraints \eqref{eq:MScons} and \eqref{eq:MSdom} are the same as in formulations $(2SP)$ and $(MSP)$. Constraints \eqref{AMSLP_d}--\eqref{AMSLP_f} capture the adaptive multistage relationship under the revision decisions through the NACs tailored for this problem for controlling the decision structures of the nodes of the scenario tree. Specifically, constraints \eqref{AMSLP_d} and \eqref{AMSLP_f} ensure the equality of state variables up to the first and after the last revision point, respectively, depending on the subtrees rooted at the revision points. Similarly, constraint \eqref{AMSLP_e} imposes equality of state variables between corresponding nodes for the remaining revision points. 

Formulation \eqref{AMSLP} is a nonlinear stochastic program  since constraints \eqref{AMSLP_d}--\eqref{AMSLP_f} depend on the set $\Gamma_i^*$ that involves revision points for every $i \in \mathcal{I}$. To linearize the formulation, we introduce an auxiliary integer variable $r_{it}$ for each $i \in  \mathcal{I}$ and $t \in \mathcal{P}$, which indicates the total number of revision points of the $i^{th}$ state variable until stage $t$. For example, in Figure~\ref{fig:mu_1_t_2}, $r_{i1}=0$ and $r_{i2}=r_{i3}=r_{i4}=1$ because the state variables are revised for the first time at stage two, and there are no revisions after that. In Figure~\ref{fig:mu_1_t_4}, $ r_{i1}=r_{i2}=r_{i3}=0$ and $r_{i4}=1$ since the first revision point is at stage four. In Figures~\ref{fig:mu_2_t_2_3} and \ref{fig:mu_2_t_2_4}, the state variables are revised twice, and the corresponding auxiliary variables are $r_{i1}=0$, $r_{i2}=1, r_{i3}=r_{i4}=2$, and $r_{i1}=0$, $r_{i2}=r_{i3}=1, r_{i4}=2$, respectively.
By leveraging the auxiliary revision variable $r_{it}$, we can reformulate the NACs 
\eqref{AMSLP_d}-\eqref{AMSLP_f} in a linear fashion. Let $\bar{x}_i$ be the upper bound on the state variable for $i \in \mathcal{I}$. 
If there is no revision between stages $t$  and $t'$ where $t' \geq t \in P$ (i.e., $r_{it} = r_{it'}$), then the NACs enforce equality of state variables at any node pair $x_{im}, x_{in}$ at stage $t'$ of the subtree $T(l)$ rooting from the node $l$ of stage $t$ for every $m, n \in S_{t^{'}} \cap \mathcal{T}(l)$. On the other hand, if there is a revision, then these constraints become inactive and serve only as an upper bound on the state variables. The generic linear formulation for AMSP with at most $\mu$ revisions is then given as follows. 
\begin{subequations}\label{AdaptiveModel}
\allowdisplaybreaks\begin{align}
(AMS_{\mu})  \ \min_{\mathbf{x}, \mathbf{y}, \mathbf{r}} \quad & \sum_{n \in \mathcal{T}} p_{n}\left( \mathbf{a_{n}}^{\top} \mathbf{x_{n}}+ \mathbf{b_{n}}^{\top} \mathbf{y_{n}}\right) \label{AdaptiveModel_obj} \\
\text{s.t.} \quad &  \eqref{eq:MScons}, \eqref{eq:MSdom} \notag \\
& x_{im} \geq x_{in}-\bar{x}_i\left(r_{it'}-r_{it} \right) \quad \forall m, n \in S_{t^{'}} \cap \mathcal{T}(l), l \in S_{t}, t^{'} \geq t, t \in \mathcal{P}, \  i \in \mathcal{I} \label{AdaptiveModel_d}\\
& x_{im} \leq x_{in}+\bar{x}_i\left(r_{it'}-r_{it} \right) \quad \forall m, n \in S_{t^{'}} \cap \mathcal{T}(l), l \in S_{t}, t^{'} \geq t, t \in \mathcal{P}, \  i \in \mathcal{I} \label{AdaptiveModel_e} \\
&  r_{i1} \leq r_{i2} \leq \dots \leq r_{iT} \leq \mu \ \ 
 \forall i \in \mathcal{I}\quad \label{AdaptiveModel_f} \\
&  r_{i(t+1)} - r_{it} \leq 1 \quad \forall t \in \mathcal{P}\backslash \{T\}, \ i \in \mathcal{I}\quad \label{AdaptiveModel_f2} \\
&  r_{i1} = 0  \quad \forall i \in \mathcal{I}\quad \label{AdaptiveModel_f3} \\
&  r_{it} \in \mathbb{Z}_{+}, \   \forall t \in \mathcal{P}, \  i \in \mathcal{I}\label{AdaptiveModel_g}
\end{align}
\end{subequations}
In this formulation, objective \eqref{AdaptiveModel_obj}, constraints \eqref{eq:MScons} and \eqref{eq:MSdom} are similar to those in formulation \eqref{AMSLP}. Constraints \eqref{AdaptiveModel_d} and \eqref{AdaptiveModel_e} ensure non-anticipativity and shape the decision structure for state variables. Constraint \eqref{AdaptiveModel_f} restricts the total number of revisions and ensures that revision variables are in non-decreasing order. Furthermore, constraints \eqref{AdaptiveModel_f2}--\eqref{AdaptiveModel_g} maintain the special structure of the integer revision variables such that there is no revision at the beginning of the planning and there can be at most one revision between consecutive stages.


We demonstrate that more revisions on the state decisions enhances the flexibility of the model, 
and we provide a formalization of this observation in the subsequent lemma.

\begin{lem}\label{Lem:binding}
For $(AMS_\mu)$, there exist an optimal solution where $r_{iT} \leq \mu$ is binding for all $i \in \mathcal{I}$.
\end{lem}
\proof{Proof:}
See Appendix \ref{App:Binding}.
\Halmos
\endproof

In the following proposition, we formalize the relationship between formulation ($AMS_\mu$) and existing optimization approaches in stochastic programming literature.

\begin{prop}\label{prop:TS-MS}
For $\mu = 0$, formulation $(AMS_\mu)$ is equivalent to two-stage stochastic programming formulation $(2SP)$. For $\mu = T - 1$, formulation $(AMS_\mu)$ is equivalent to multistage stochastic programming formulation $(MSP)$. 
\end{prop}

\proof{Proof:} See Appendix \ref{App:TS-MS}. \Halmos 
\endproof 

Proposition \ref{prop:TS-MS} highlights that formulation $(AMS_\mu)$ is a generalization of existing stochastic programming approaches. Furthermore, when $\mu$ is set to 1,  formulation reduces to the adaptive two-stage stochastic programming proposed in \cite{basciftci2019adaptive}. 

By increasing $\mu$, we can enhance the flexibility in decision-making, leading to improvements in the objective value. We formalize this in the following proposition, where we let $z(\cdot)$ denote the optimal objective value of a given problem. For example, $z(AMS_\mu)$ is the optimal objective value of formulation $(AMS_\mu)$. 
It is important to note that this study assumes that formulation $(2SP)$ is feasible and consequently $(AMS_\mu)$ is feasible under any $\mu$ value.

\begin{prop} \label{prop:ObjectiveComparison}
$z(2SP) =  z(AMS_0) \geq z(AMS_1) \geq \dots \geq z(AMS_{(T-2)}) \geq z(AMS_{(T-1)}) = z(MSP)$. 
\end{prop}

\proof{Proof:} See Appendix \ref{App:ObjectiveComparison}.\Halmos
\endproof

We next analyze the complexity of this problem and formalize it in the following proposition. 
\begin{prop}\label{prop:NPHard}
AMSP 
is NP-Hard. 
\end{prop}

\proof{Proof:} 
See Appendix \ref{App:NPHard}. 
\Halmos
\endproof

A special case of AMSP arises when the state variables are binary. This characteristic permits the relaxation of the integrality requirement for the revision decision variables $\mathbf{r}$. We formalize this observation in Proposition~\ref{prop:binary}.  

\begin{prop} \label{prop:binary} Let the state variables be binary with $x_{in} \in \{0, 1\}$ for all $i \in \mathcal{I}, \> n \in \mathcal{T}$. Then, the integrality conditions on the revision decision variables $r_{it}$ for all $i \in \mathcal{I}$ and $t \in \mathcal{P}$ can be relaxed within $(AMS_\mu)$.
\end{prop} 

\proof{Proof:}
See Appendix \ref{App:binary}.
\Halmos
\endproof

\section{Elimination Techniques for Non-anticipaticity Constraints}\label{sec:NACsElimination}

The size of formulation $(AMS_\mu)$ increases exponentially with the increasing number of stages, making the model challenging to solve. In this section, we enhance formulation $(AMS_\mu)$ by identifying a series of constraint elimination techniques for the NACs 
\eqref{AdaptiveModel_d} and \eqref{AdaptiveModel_e}. Recall that these constraints consider any combination of stages $t' \geq t \in \mathcal{P}$, and impose equality between nodes $m, n$ at stage $t'$ if there exists a subtree rooted at time $t$ containing both nodes $m$ and $n$ when there is no revision between $t$ and $t'$.  
In the following, we first present a series of propositions that identify how we can reduce the number of NACs, and then combine these results to effectively reformulate $(AMS_\mu)$. 

The first elimination, presented in Proposition \ref{Prop:Elimination1}, results from imposing a cyclic relation between the nodes associated with constraints \eqref{AdaptiveModel_d} and \eqref{AdaptiveModel_e}. Specifically, instead of considering each $m, n$ pair in the subtree sets $S_{t^{'}} \cap \mathcal{T}(l)$ for every $l \in S_t$, we order the nodes in these sets for each $t'$, and enforce equality between these ordered pair of nodes.
\begin{prop}\label{Prop:Elimination1}
Let $t \in \mathcal{P}$, $t^{'} \geq t$, and $l \in S_{t}$ be given. Let $K = |S_{t^{'}} \cap \mathcal{T}(l)|$. 
Without loss of generality, we can order the nodes in $S_{t^{'}} \cap \mathcal{T} (l)$ as $(n, n +1, \dots,n + K-1)$. 
Then, it is sufficient to include constraints \eqref{AdaptiveModel_d} and \eqref{AdaptiveModel_e} only for the consecutive node pairs (e.g., $(n,n+1),(n+1,n+2),\dots,(n+K-1,n)$) from the set $S_{t^{'}} \cap \mathcal{T}(l)$ as follows: 
\begin{align}\label{eq:Property1}
\begin{split}
 x_{in} \geq x_{i(n+1)} -\bar{x}_i\left(r_{it'}-r_{it} \right) \quad \quad n \in  S_{t^{'}} \cap \mathcal{T}(l),\quad l \in S_t,\quad t'\geq t,\quad t \in \mathcal{P},\quad i \in \mathcal{I} \\
\end{split}
\end{align}
We note that we should remove the ordered node $n + K$ from the set $S_{t^{'}} \cap \mathcal{T} (l)$ and add the constraint between $x_{i(n+K-1)}$ and $x_{in}$, but to enhance readability we slightly abuse notation in constraints \eqref{eq:Property1} by using node $n + K$ for node $n$.
\end{prop}

\proof{Proof:}
See Appendix \ref{App:Elimination1}.\Halmos
\endproof
In Proposition \ref{Prop:Elimination1}, we focus on node pairs $(m,n)$ that belong to the same stage $t'$, and remove redundant constraints by exploiting an arbitrary ordering argument between such nodes. Alternatively, we observe that any given node pair $(m, n)$ belonging to the set $S_{t'}$ may be included in various NACs of different $(t,t')$ combinations. However, it is sufficient to link a node pair $(m, n)$ depending on only one particular $(t,t')$ combination, where $t$ is the stage of the last common ancestor of nodes $(m,n)$. This observation is based on the non-decreasing nature of the $r$ variables as formalized below. 

\begin{prop}\label{Prop:Elimination2}
Let $t_{a(m,n)}$ be the stage of the last common ancestor of node pair $(m,n)$ from the set $S_t$ of stage $t$. It is sufficient to include NACs of $(m,n)$ only for $t, t_{a(m,n)}$ combination, and accordingly, constraints \eqref{AdaptiveModel_d} and \eqref{AdaptiveModel_e}  can be replaced by the following constraints:
\begin{subequations} \label{C:elimination2}
\begin{align}
& x_{im} \geq x_{in}-\bar{x}_i\left(r_{it}-r_{it_{a(m,n)}} \right) \quad \forall m, n \in S_{t},  t \in \mathcal{P}, i \in \mathcal{I} \\
& x_{im} \leq x_{in}+\bar{x}_i\left(r_{it}-r_{it_{a(m,n)}} \right) \quad \forall m, n \in S_{t},   t \in \mathcal{P},  i \in \mathcal{I}  
\end{align}\end{subequations}
\end{prop}

\proof{Proof:}
See Appendix \ref{App:Elimination2}.
\Halmos 
\endproof

In Proposition \ref{Prop:Elimination2}, we include only one pair of NACs for any $(m, n)$ node pair in the same stage. Furthermore, we observe that NACs are not required for certain $(m, n)$ pairs, depending on the value of $\mu$. This observation is based on the fact that the model revises state decisions to increase flexibility as much as possible, resulting in $r_{iT}^* = \mu$ for $i \in \mathcal{I}$ at an optimal solution, as shown in Lemma \ref{Lem:binding}.

\begin{prop} \label{Prop:Elimination3}
Let $t\in \{1,\dots,T\}$ and $t'\geq t$. For given $\mu$, NACs between a particular $(t',t)$ combination become redundant, if $t'>t+(T-\mu)-1$. Then, constraints \eqref{AdaptiveModel_d} and \eqref{AdaptiveModel_e} reduce to:
\begin{subequations}
\allowdisplaybreaks\begin{align}
x_{i m} \geq x_{i n}-\bar{x}_i\left(r_{it'}-r_{it} \right) \ \forall m, n \in S_{t^{'}} \cap \mathcal{T}(l), l \in S_{t},  \min(t+ (T-\mu)-1, T) \geq t^{'} \geq t, t \in \mathcal{P}, i \in \mathcal{I} \label{eq:El31}\\
x_{i m} \leq x_{i n}+\bar{x}_i\left(r_{it'}-r_{it} \right) \  \forall m, n \in S_{t^{'}} \cap \mathcal{T}(l), l \in S_{t}, \min(t+ (T-\mu)-1, T) \geq t^{'} \geq t, t \in \mathcal{P}, i \in \mathcal{I} \label{eq:El32}
\end{align}
\end{subequations}
\end{prop}

\proof{Proof:}
See Appendix \ref{App:Elimination3}. \Halmos
\endproof

We remark that Propositions~\ref{Prop:Elimination1}-\ref{Prop:Elimination3} do not dominate each other and thus the number of NACs reduces further if we combine the propositions as follow:
\begin{equation}
 x_{in} \geq x_{i(n+1)}-\bar{x}_i\left(r_{it}-r_{it_{a(n,n+1)}} \right) \quad \forall n \in S_{t}, t_{a(n,n+1)} \in \mathcal{P}_t,  t \in \mathcal{P}, i \in \mathcal{I} \label{combinedNACs}  
\end{equation}
where set $P_t$ is defined as $\{(t-(T-\mu)+1)^+,\dots, t-1\}$ for a given $t$.  We summarize the main outcome of this section in Theorem \ref{th:theoneandonlytheorem}. 

\begin{theorem} \label{th:theoneandonlytheorem}
AMSP can be correctly formulated as follows: 
\begin{subequations}\label{ReducedAdaptiveModel}
\allowdisplaybreaks\begin{align}
\quad \quad \quad \min_{\mathbf{x}, \mathbf{y}, \mathbf{r}} \quad & \sum_{n \in \mathcal{T}} p_{n}\left( \mathbf{a_{n}}^{\top} \mathbf{x_{n}}+ \mathbf{b_{n}}^{\top} \mathbf{y_{n}}\right) \label{RAdaptiveModel_obj} \\
\text{s.t.} \quad &  \eqref{eq:MScons}, \eqref{eq:MSdom}, \eqref{AdaptiveModel_f}-\eqref{AdaptiveModel_g}, \eqref{combinedNACs}  \notag  
\end{align}
\end{subequations}
We refer to this formulation as $(AMS_{\mu}+)$, which is equivalent to formulation $(AMS_\mu)$.
\end{theorem}

We demonstrate the impact of Propositions \ref{Prop:Elimination1}--\ref{Prop:Elimination3} on the number of necessary NACs by considering an example scenario tree with ten stages $(T = 10)$ and two branches $(B = 2)$ focusing on a single state variable ($I=1$). Table \ref{tab:Normal NA} displays the required number of NACs for each $(t, t')$ pair in formulation $(AMS_\mu)$ which does not include the proposed elimination techniques. We observed a total of 688810 NACs, which has notable implications for the complexity of the model.

Table~\ref{tab:TNACs} displays the resulting changes in the number of included NACs when we apply Propositions \ref{Prop:Elimination1}-\ref{Prop:Elimination3}. Applying Proposition~\ref{Prop:Elimination1} alone reduces the total number of NACs to 8194 by imposing a cyclic relationship on the nodes rather than having NACs for every possible combination of $m,n$ pairs in stage $t'$. Consequently, the number of NACs for every $t'$ equals the number of nodes in that stage. By applying both Propositions \ref{Prop:Elimination1} and \ref{Prop:Elimination2}, we further reduce the total number of NACs to 2026. Prior to this, a node pair could be seen in many constraints, but now it is unique to a particular $t,t'$ pair. Depending on the maximum allowed number of revisions $\mu$, we can further decrease the number of constraints by considering Proposition~\ref{Prop:Elimination3}. We display an example where $\mu$ is set to 4. In this case, 52 more constraints ca be eliminated. The number of eliminated constraints changes to 2, 8, 22, 114, 240, 494, and 1004 when $\mu$ is set to be 1, 2, 3, 5, 6, 7, and 8, respectively.

\begin{table}[h]
\small
        \centering
        \caption{Number of NACs for formulation $(AMS_\mu)$ with $B=2$ Branches } \label{tab:Normal NA}
        \begin{tabular}{l|rrrrrrrrr}
         \toprule \backslashbox{$t$}{$t'$}  &2 &3 &4 &5 & 6 & 7 & 8 & 9 & 10 \\
        \midrule
       1 &2 &12 &56 &240 &992 &4032 &1.6e4  &6.5e4  &2.6e5 \\
       2 &  &4  &24 &112 &480 &1984 &8064  &3.2e4  &1.3e5 \\
       3 &  &   &8  &48  &224 &960  &3968  & 1.6e4 &6.5e4 \\
       4 &  &   &   &16  &96 &448   &1920  &7936  & 3.2e4 \\
       5 &  &   &   &    &32	&192   &896	  &3840	 &1.6e4  \\
       6 &  &   &   &    &   &64    &384   &1792  &7680 \\
       7 &  &   &   &    &   &      &128   &768   &3584 \\
       8 &  &   &   &    &   &      &      & 256   & 1536 \\
       9 &  &   &   &    &   &      &    &        & 512 \\ \bottomrule 
\end{tabular} 
\end{table}

\begin{table}[h]
        \centering
           \caption{Number of NACs with Propositions \ref{Prop:Elimination1} - \ref{Prop:Elimination3}}
           \label{tab:TNACs}
        \resizebox{\columnwidth}{!}{%
       \begin{tabular}{l|rrrrrrrrr|rrrrrrrrr|rrrrrrrrr} \toprule 
        & \multicolumn{9}{c}{Proposition~\ref{Prop:Elimination1}} & \multicolumn{9}{|c|}{Propositions~\ref{Prop:Elimination1} and \ref{Prop:Elimination2}} & \multicolumn{9}{c}{ Propositions~\ref{Prop:Elimination1}, \ref{Prop:Elimination2},  and \ref{Prop:Elimination3} ($(AMS_\mu+)$)} \\
        \midrule 
         \backslashbox{$t$}{$t'$}  &2 &3 &4 &5 & 6 & 7 & 8 & 9 & 10 &2 &3 &4 &5 & 6 & 7 & 8 & 9 & 10 &2 &3 &4 &5 & 6 & 7 & 8 & 9 & 10 \\
        \midrule 
      1 &2 &4 &8 &16 &32&64 &128 &256 &512  &2 &2 &2 &2  &2  &2  &2   &2   &2  &2 &2 &2 &2  &2  &  &   &   & \\
       2 &  &4 &8 &16 &32&64 &128 &256 &512 &  &4 &4 &4  &4  &4  &4   &4   &4 &  &4 &4 &4  &4  &4  &   &   & \\
       3 &  &  &8 &16 &32&64 &128 &256 &512 &  &  &8 &8  &8  &8  &8   &8   &8 &  &  &8 &8  &8  &8  &8   &   &\\
       4 &  &  &  &16 &32&64 &128 &256 &512 &  &  &  &16 &16 &16 &16  &16  &16 &  &  &  &16 &16 &16 &16  &16  & \\
       5 &  &  &  &   &32&64 &128 &256 &512 &  &  &  &   &32 &32 &32  &32  &32 &  &  &  &   &32 &32 &32  &32  &32\\
       6 &  &  &  &   &  &64 &128 &256 &512 &  &  &  &   &   &64 &64  &64  &64 &  &  &  &   &   &64 &64  &64  &64\\
       7 &  &  &  &   &  &   &128 &256 &512 &  &  &  &   &   &   &128 &128 &128 &  &  &  &   &   &   &128 &128 &128\\
       8 &  &  &  &   &  &   &    &256 &512 &  &  &  &   &   &   &    &256 &256 &  &  &  &   &   &   &    &256 &256 \\
       9 &  &  &  &   &  &   &    &    &512 &  &  &  &   &   &   &    &    &512 &  &  &  &   &   &   &    &    &512\\ \bottomrule
       \end{tabular}}

\end{table}
       
The total number of NACs for formulation $(AMS_\mu)$ can be calculated as $I \times \sum_{t=1}^{T-1} \sum_{j=1}^{T-t} B^{t-1} \binom{B^j}{2} \times 2$. However, after applying all of the propositions, the total number of NACs in formulation $(AMS_\mu+)$ reduces $I \times \sum_{t=1}^{T-1} \sum_{j=1}^{\min(T-t,t+T-\mu-1)} B^{t}$, which represents a substantial decrease in the number of NACs. Table \ref{tab:Number of NACs} illustrates the total number of NACs for varying numbers of stages and branches when there is a single state variable. Note that Proposition \ref{Prop:Elimination3} depends on the value of $\mu$, and we present two examples when $\mu$ is set to 2 and 4. The propositions lead to a significant reduction in the total number of NACs, ranging from 90\% to almost 100\%. Note that the total number of NACs in Table \ref{tab:Number of NACs} for the case of 2 branches, 10 stages, and $\mu=4$ corresponds to the sum of the NACs listed in Table \ref{tab:TNACs} for $(AMS_\mu+)$ when all propositions are applied.

\vspace{-2mm}
\begin{table}[h]
\caption{Total number of NACs over different instance sizes} \label{tab:Number of NACs}
\centering
\resizebox{\columnwidth}{!}{%
\begin{tabular}{l| rrrrrr|rrrrrr} \toprule
Branch & \multicolumn{6}{c|}{2} & \multicolumn{6}{c}{3}\\ \midrule
Stage & 5 & 6 & 7 & 8 & 9 & 10 & 5 & 6 & 7 & 8 & 9 & 10 \\ \hline
Formulation $(AMS_\mu)$ & 522 & 2346 & 1.0e4 & 4.2e4 & 1.7e5 & 6.9e5 & 1.0e4 & 9.7e4 & 8.9e5 & 8.0e6 & 7.3e7 & 6.5e8 \\
Formulation $(AMS_\mu+)$, $\mu=2$ & 44 & 106 & 232 & 486 & 994 & 2018 & 159 & 522 & 1614 & 4893 & 1.5e4 & 4.4e4 \\
Formulation $(AMS_\mu+)$, $\mu=4$ & 0 & 62 & 188 & 442 & 952 & 1974 & 0 & 363 & 1455 & 4734 & 1.4e4 & 4.4e4 \\ \bottomrule
\end{tabular}}
\end{table}

\vspace{-6mm}
\section{Solution Approach}
\label{sec:SolutionApproach}
This section introduces a general solution approach for solving AMSPs with mixed-integer state and stage variables. Standard solution approaches for multistage stochastic programming cannot be used due to the presence of integer revision decisions $r_{it}$ for every $i \in \mathcal{I}$ and $t \in \mathcal{P}$. To address this issue, we propose a solution approach that separates these complicating decision variables from the remainder of the problem, which leverages Benders decomposition \citep{benders} and the integer L-shaped method \citep{van1969shaped}. More specifically, we decompose the problem into a master problem, in which we make the revision decisions, and a subproblem, in which we solve the adaptive multistage problem under the given revision decisions. Then, we present an iterative process, where we first solve the master problem and then add cuts to this problem using the subproblem. Note that integer L-shaped cuts are designed for pure binary master programs, while Benders cuts require convexity of the subproblem. In order to accommodate the presence of integer revision decisions in the master problem and possible mixed-integer subproblems, we make the necessary modifications to incorporate these cuts, as explained later in this section.

First, we present the master problem \eqref{MasterProblem} by introducing an auxiliary variable $\theta$ to approximate the cost of the subproblem as follows:
\vspace{-3mm}
\begin{subequations} \label{MasterProblem}
\allowdisplaybreaks\begin{align}
(MP) \ \min_{\mathbf{r}, \theta} \quad &  \ \theta \\
\text{s.t.} \quad & \ r_{i1} \leq r_{i2} \leq \dots \leq r_{iT} \leq \mu  &\forall& i \in \mathcal{I} \label{Master:Adaptive1} \\
&  r_{i(t+1)} - r_{it} \leq 1  &\forall& t \in \mathcal{P}\backslash \{T\}, \ i \in \mathcal{I}\quad \label{Master:Adaptive12} \\
&  r_{i1} = 0   &\forall& i \in \mathcal{I}\quad \label{Master:Adaptive13} \\
&  r_{it} \in \mathbb{Z}_{+}  &\forall& t \in \mathcal{P},  i \in \mathcal{I} \label{Master:Adaptive2} \\
&  \theta \geq  L \label{Master:LB} \\
&  (r, \theta) \in \Phi \label{Master:Optimality} \\
&  (r, \theta) \in \Psi \label{Master:BendersOptimality} 
\end{align}
\end{subequations}

\vspace{-2mm}
Here, constraints \eqref{Master:Adaptive1}-\eqref{Master:Adaptive2} correspond to the integer revision point related constraints \eqref{AdaptiveModel_f}-\eqref{AdaptiveModel_g} as in models \eqref{AdaptiveModel} and \eqref{ReducedAdaptiveModel}. Constraint \eqref{Master:LB} sets an initial lower bound on the variable $\theta$. To set this lower bound $L$, we solve the multistage problem as a preprocessing step since the objective of the multistage stochastic program provides a lower bound on the adaptive multistage problem. Sets $\Phi$ and $\Psi$ in constraints \eqref{Master:Optimality} and \eqref{Master:BendersOptimality} correspond to the sets of L-shaped and Benders optimality cuts, respectively. We incorporate only optimality cuts and omit feasibility cuts, given our assumption that the AMSP is feasible under any revision vector $\mathbf{r}$ (see Proposition \ref{prop:ObjectiveComparison}).

We let $\mathbf{\bar r}$ denote the vector of revision decisions obtained from the master problem. Then, 
using the improved formulation of the adaptive multistage problem $(AMS_\mu+)$, the subproblem can be defined as follows:
\vspace{-2mm}
\begin{subequations}\label{SubProblem}
\allowdisplaybreaks\begin{align}
(SP) \ \mathcal{Q}(\mathbf{\bar r}) = \min_{\mathbf{x}, \mathbf{y}} \quad & \sum_{n \in \mathcal{T}} p_{n}\left(\mathbf{a_{n}}^{\top} \mathbf{x_{n}}+\mathbf{b_{n}}^{\top} \mathbf{y_{n}}\right) \label{SP:obj} \\
\text{s.t.} \quad &  \sum_{m \in P(n)} \mathbf{C_{n m}} \mathbf{x_{m}}+\mathbf{D_{n}} \mathbf{y_{n}} \geq \mathbf{d_{n}}  &\forall& n \in \mathcal{T} \label{SP:c} \\ 
&  \mathbf{x_{n}} \in \mathcal{X}_{n}, \mathbf{y_{n}} \in \mathcal{Y}_{n}, \quad &\forall& n \in \mathcal{T} \label{SP:set}  \\
&  x_{in} \geq x_{i(n+1)}-\bar{x}_i\left(\bar r_{it}-\bar r_{it_{a(n,n+1)}} \right) &\forall& n \in S_{t}, t_{a(n,n+1)} \in \mathcal{P}_t,  t \in \mathcal{P}, i \in \mathcal{I} \label{SP:r1} 
\end{align}
\end{subequations}

The integer L-shaped method can be applied to problems where the master problem includes pure binary decision variables \citep{laporte1993integer}. However, our master problem includes integer revision decisions. One potential solution is representing these integer variables as a sum of binary decision variables. This representation allows expressing the optimality cuts obtained from the subproblem in terms of these auxiliary binary variables, which introduces additional complexity to the solution procedure. Instead, we propose an alternative scheme by utilizing the special structure of the revision decisions. Specifically, the difference in revision points between two consecutive stages, denoted as $(r_{it}-r_{i(t-1)})$ for $t \in \mathcal{P} \setminus \{1\}$ and $i \in \mathcal{I}$, can only be one or zero. If this value is one, then there is a revision at time $t$ for the state variable $i$, and zero, otherwise. We utilize this information to define the set $Y_i\left(\mathbf{\bar r} \right):=\left\{t \in \mathcal{P}: \bar r_{it}-\bar r_{i(t-1)}=1\right\}$ under a given revision vector $\mathbf{\bar r}$, representing the set of periods with revision for the state variable $i$. Accordingly, using the subproblem \eqref{SubProblem}, the integer L-shaped based optimality cuts can be obtained as follows: 
\begin{equation}\label{LshapedCut}
    \theta \geq\left(\mathcal{Q}\left(\mathbf{\bar r}\right)-L\right) \left(\sum_{i \in \mathcal{I}} \sum_{t \in Y_i\left(\bar r \right)}\left(r_{it}-r_{i(t-1)}-1\right)- \sum_{i \in \mathcal{I}} \sum_{t \notin Y_i\left(\bar r \right)} (r_{it}-r_{i(t-1)})\right)+\mathcal{Q}\left(\mathbf{\bar r}\right)
\end{equation}

Since set $\Phi$ contains exponentially many cuts in the form of \eqref{LshapedCut}, obtaining these cuts requires solving mixed-integer subproblems, which may cause a computational burden. Therefore, we integrate Benders optimality cuts to decrease the feasible region of the master problem. As convexity of the subproblems is necessary for Benders decomposition to converge to an optimal solution and derive optimality cuts  \citep{geoffrion1972generalized},  we convexify the subproblem by relaxing the integrality of the state and stage variables. We refer to the relaxed version of formulation $(SP)$ as $(RSP)$ and denote its objective function as $\underline{Q}(\bar r)$. Accordingly, by taking the dual of this relaxed subproblem, we obtain Benders optimality cuts $\Psi$ as follows:
\begin{equation} \label{eq:BendersOptCut}
\theta \geq \sum_{n \in \mathcal{T}}  \mathbf{d_{n}}^{\top}  \mathbf{\bar{\Lambda}_n} -  \sum_{i \in \mathcal{I}} \sum_{t \in \mathcal{P}} \sum_{n \in \mathcal{S}_t} \sum_{t_{a(n,n+1)} \in \mathcal{P}_t}  \bar{x}_i\left(r_{it}-r_{it_{a(n,n+1)}} \right) \bar{\Pi}_{itn(n+1)}
\end{equation}
where $\bar{\mathbf{\Lambda}}$ and $\bar{\mathbf{\Pi}}$ represents the optimal values of the dual variables corresponding to constraints \eqref{SP:c} and \eqref{SP:r1}, respectively, under the solution $\bar r$. For simplicity in notation, we assume that constraints \eqref{SP:set} are non-negativity restrictions, and we omit their corresponding dual variables in illustrating \eqref{eq:BendersOptCut}. Solving the dual of $(RSP)$ and obtaining Benders cuts is computationally easier than solving $(SP)$ and obtaining L-shaped cuts, as it involves linear programs rather than their mixed-integer counterparts. Therefore, we solve $(RSP)$ and incorporate Benders cuts in each iteration as a revision vector is obtained by the master problem. However, the decision to solve $(SP)$ and include L-shaped cuts depends on the improvement observed in the solution of $(RSP)$ using the obtained revision vector. In this way, we selectively add L-shaped cuts in a subset of the iterations for the sake of solution time. The detailed procedure is further explained later in Algorithm \ref{SolutionAlgorithm}.

To accelerate the solution algorithm further, we examine the timing of the revision point decisions over the problem instances with different lengths of planning horizons.

\begin{obs} \label{obs:Preprocessing}
In our preliminary experiments, we 
consider a modified version of the original problem $(AMS_\mu)$, by reducing the planning horizon by $n$ stages, resulting in a shortened scenario tree denoted as $\mathcal{T}(n_0,T-n)$, where $n_0$ is the root node of $\mathcal{T}$. Specifically, we remove the stages and decision variables from $T-n+1$ to $T$ from the original problem, and solve the reduced version, which we call $(AMS_\mu^{T-n})$.
In comparison to the solution obtained from the original problem $(AMS_\mu)$, we have observed that as the planning horizon is reduced, revision points of the original problem that are scheduled between stages $1, \dots, T-n$ do not shift to subsequent stages in the solution obtained from $(AMS_\mu^{T-n})$; instead, they either remain stable or shift to the previous stages.
\end{obs}

We leverage this observation in the preprocessing step of our algorithm to eliminate a set of candidate feasible solutions by solving a simpler variant of the original problem. 
Therefore,  we first solve $AMS_{\mu}^{T-n}$ with a reduced planning horizon for a given $n \in \{1,\dots, T-1\}$. We then obtain an optimal revision solution from $AMS_{\mu}^{T-n}$, denoted as $\mathbf{\underline r}$. Due to Observation \ref{obs:Preprocessing}, we expect the revision points of the original problem at stages $1, \dots, T-n$ to either remain at the same stage or move to later stages compared to the shorter-term version of this problem. Thus, we propose the heuristic cuts \eqref{Heuristic_Cuts} to be added to the master problem. 
\begin{equation} \label{Heuristic_Cuts}
    r_{it} \geq \underline r_{it} \quad \forall t \in \{1,\dots,T-n\}, i \in \mathcal{I}
\end{equation}
    
\begin{algorithm}[h]
\small
\caption{Solution Algorithm}
\label{SolutionAlgorithm}
\begin{algorithmic}[1]
\STATE Set LB = -$\infty$, UB = RUB =$\infty$, $\Psi = \emptyset$, $\Phi = \emptyset$.
\STATE  $L \leftarrow$ Solve $(MSP)$ to obtain $z(MSP)$.
\STATE Solve $(AMS_{\mu}^{T-n}+)$ to obtain $\mathbf{\underline r}$.
\STATE Add heuristic cuts \eqref{Heuristic_Cuts} to $(MP)$.  \label{Algortihm_Heuristic}
\STATE Set $\mathbf{\bar r}$ to $\mathbf{\underline r}$.
\WHILE {$(1-\text{LB}/\text{UB})<\epsilon $}
\STATE  $z(RSP) \leftarrow$ Solve $(RSP)$ subject to $\mathbf{\bar r}$.
\STATE Generate Benders optimality cuts \eqref{eq:BendersOptCut} and add them to $\Psi$.
\IF{$z(RSP) \leq $ RUB} \label{Algortihm_RUB}
\STATE RUB $\leftarrow z(RSP)$.
\STATE  $z(SP) \leftarrow$ Solve $(SP)$ subject to $\mathbf{\bar r}$.
\STATE Generate L-shaped optimality cuts \eqref{LshapedCut} and add them to $\Phi$.
\IF{$z(SP) \leq $ UB}
\STATE UB $\leftarrow z(SP)$.
\STATE $\mathbf{r^*} \leftarrow \mathbf{\bar r}$.
\ENDIF
\ENDIF
\STATE Solve $(MP)$ with $\Psi$ and $\Phi$ to obtain $\mathbf{\bar r}$.
\STATE Set LB $\leftarrow z(MP)$.
\ENDWHILE 
\end{algorithmic}
\end{algorithm}
\normalsize

Combining above, we summarize our solution approach in Algorithm \ref{SolutionAlgorithm}. In the first step of this algorithm, we solve the multistage stochastic problem to use its objective as a lower bound for $(MP)$ within constraint \eqref{Master:LB}. Next, we solve the $(AMS_{\mu}+)$ problem for $T-n$ stages to obtain heuristic cuts and an initial solution $\mathbf{\underline r}$. We then include heuristic cuts to the master problem. In the first step of the while loop, we solve the relaxed subproblem under fixed $\mathbf{\underline r}$ to obtain its objective value $z(RSP)$ and generate Benders optimality cuts to be added to $\Psi$. If the value of $z(RSP)$ is less than or equal to relaxed upper bound (RUB), the upper bound (UB) is updated to the value of $z(RSP)$. The algorithm then solves the subproblem under fixed  $\mathbf{\underline r}$ and generates L-shaped optimality cuts to be added to $\Phi$. If the value of $z(SP)$ is less than or equal to UB, the upper bound is updated to the value of $z(SP)$, and $\mathbf{r^*}$ is set to  $\mathbf{\underline r}$. We then solve $(MP)$ with the obtained cuts $\Psi$ and $\Phi$ to get revisions  $\mathbf{\underline r}$ and set the lower bound (LB) to $z(MP)$. We repeat this procedure until the gap (1-LB/UB) is smaller than $\epsilon$.

\begin{remark}
If the heuristic cut in Step \ref{Algortihm_Heuristic} and the condition in Step \ref{Algortihm_RUB} of Algorithm \ref{SolutionAlgorithm} are excluded, the proposed approach is guaranteed to converge to an optimal solution in a finite number of steps. However, our computational experiments have demonstrated that incorporating these two steps significantly accelerates the algorithm's convergence and results in an optimal solution for almost all instances.
\end{remark}

\section{Computational Experiments}\label{sec:ComputationalExperiments}
We demonstrate the value of AMSP by addressing two distinct multistage stochastic optimization problems that require a balance between decision flexibility and decision commitment. First, in Section \ref{sec:LotSizing}, we study a canonical optimization problem from supply chain management; the stochastic uncapacitated lot-sizing problem. Second, in Section \ref{Sec:GenerationExpansion}, we provide a case study of the generation expansion planning problem, which is fundamental in our society's quest toward a transition to renewable energy under various uncertainties. Besides showing the value of AMSP, we demonstrate the computational performance of the NACs elimination techniques proposed in Section \ref{sec:NACsElimination} and our solution algorithm presented in Section \ref{sec:SolutionApproach}. 

All the computational experiments are carried out on an Intel Xeon E5 2680v3 CPU with a 2.5 GHz processor and 32 GB RAM. The implementation is performed in Python, utilizing the \texttt{Gurobi} 9.1.0 solver to solve the mathematical models \eqref{AdaptiveModel} and \eqref{ReducedAdaptiveModel}, as well as the master problem \eqref{MasterProblem} and subproblem \eqref{SubProblem} within our solution algorithm. For consistency, a time limit of 5 hours is imposed on all experiments.

To evaluate the benefit of AMSP, we introduce a metric called the ``Value of AMSP" (VAMS), which compares its solution quality with two-stage $(2SP)$ and multistage stochastic programming $(MSP)$ solutions. Note that two-stage stochastic programming $(2SP)$ offers a feasible solution for AMSP and provides an upper bound, while multistage stochastic programming $(MSP)$ is a relaxation for AMSP, providing a lower bound. Accordingly, we define the VAMS for a given level of flexibility $\mu$ as follows:
\begin{equation}
    \text{VAMS}=\frac{z(2SP)-z(AMS_\mu)}{z(2SP)-z(MSP)} \times 100 \%
\end{equation}

The VAMS metric quantifies the proximity between the objective values of AMSP and multistage stochastic programming relative to two-stage stochastic programming. Notably, the VAMS value is 0\% if $\mu=0$, as $z(AMS_{0})$ is equal to the optimal value $z(2SP)$ of two-stage stochastic programming. Conversely, the VAMS value is 100\% if $\mu=T-1$, as $z(AMS_{T-1})$ is equal to the optimal value $z(MSP)$ of multistage stochastic programming. A higher VAMS implies that AMSP can generate solutions comparable to those of multistage stochastic programming by optimizing revision points.



\subsection{Adaptive Multistage Stochastic Lot-Sizing Problem}\label{sec:LotSizing}
The stochastic lot-sizing problem \citep{guan2008polynomial} concerns determining optimal production schedules of a given number of production sources to cover the demand of a single item over a finite planning horizon. The goal is to minimize the expected sum of production setup, production amount, and inventory holding cost considering stochasticity of demand and cost parameters. The production setup decisions are made first to determine whether to initiate production for each source at a given stage. Based on the setup decisions, the production amount for each source and the associated inventory levels are determined. Multistage stochastic programming is commonly used for lot-sizing problems; it provides flexibility to update all decisions based on the realized uncertainty through the planning horizon. However, in many practical settings, frequent updates to decisions can be detrimental to production scheduling performance, leading to decreased productivity and increased global costs \citep{herrera2016reactive}. Furthermore, such frequent updates might be even impractical due to various factors, such as dependencies on predetermined schedules for other items \citep{kouvelis2006supply}. Addressing such instabilities while enabling flexible production management is a focus of various studies  \citep{herrera2016reactive,meistering2017stabilized}.
Thus, optimizing the stages at which the  decisions will be revised has a significant potential to enhance this problem setting. Therefore, in this section, we introduce the Adaptive Multistage Stochastic Uncapacitated Lot-Sizing Problem. 

We consider a set of sources represented by $i \in \mathcal{I}$ and a scenario tree $\mathcal{T}$ that describes the underlying stochastic process. Each scenario in the tree has an equal probability of occurring. The system's uncertainty lies in setup costs $\alpha_{in}$, production costs $\beta_{in}$, inventory holding costs $h_n$, and demand $d_n$ for every node $n \in \mathcal{T}$. 
We model the setup decisions as binary variables $x_{in}$, which take a value of one if there is a setup for production of source $i$ at node $n$, and zero otherwise. We represent the production amount of source $i$ at node $n$  and the inventory level at node $n$ as continuous variable $y_{in}$ and $s_n$, respectively. We denote $M$ as a large number for constraint linearization purposes. Note that the setup variables are the \textit{state} decisions, while production and inventory holding variables are the \textit{stage} decisions. Based on the multistage stochastic, single-item, uncapacitated lot-sizing problem formulation from \citep{guan2008polynomial}, we propose the following AMSP variant:
\begin{subequations}\label{model:MLS}
\allowdisplaybreaks
\begin{align}
\min \quad &  \sum_{n \in \mathcal{T}} p_{n}   \left( \sum_{i \in \mathcal{I}} \left(\alpha_{in} x_{in}+\beta_{in} y_{in} \right)+ h_n s_n \right)  \label{MLS_obj} \\
\text{s.t.} \quad & s_{a(n)}+\sum_{i \in \mathcal{I}} y_{in}=d_n+s_n  &\forall& n \in \mathcal{T} \label{MLS_stock}  \\
&  y_{in} \leq M x_{in}   &\forall& n \in \mathcal{T}, i \in \mathcal{I} \label{MLS_setup} \\
& x_{in} \in\{0,1\},  y_{in}, s_n \in \mathbb{R}_{+}    & \forall& n \in  \mathcal{T}, i \in \mathcal{I}  
\label{MLS_set} \\
&\eqref{AdaptiveModel_f} - \eqref{AdaptiveModel_g}, \eqref{combinedNACs} \notag
\end{align}
\end{subequations}
Here, objective \eqref{MLS_obj} minimizes the total expected cost of production, setup, and inventory holding. Constraint \eqref{MLS_stock} ensures the balance of the inventory flow. Constraint \eqref{MLS_setup} controls the setup-production relationship. Constraint \eqref{MLS_set} defines feasibility sets of the variables. Finally, constraints \eqref{AdaptiveModel_f} - \eqref{AdaptiveModel_g}, \eqref{combinedNACs} are the necessary NACs, which also include the revision points $\mathbf{r}$.

To analyze the stochastic lot-sizing problem described above, we create a scenario tree by generating random cost and demand parameters at each node. Specifically, the setup cost $\alpha_{in}$ is drawn from $\mathcal{U}(100,250)$, production cost $\beta_{in}$ is drawn from $\mathcal{U}(0,5)$, inventory holding cost $h_n$ is drawn from $\mathcal{U}(0,5)$, and demand $d_n$ is drawn from $\mathcal{U}(0,100\times I)$. We create ten different scenario trees and report the results on average.

\subsubsection{Potential of AMSP.}
\begin{figure}[h]
\centering
    \caption{The performance of AMSP}
    \label{fig:Effect of mu}
\begin{subfigure}{.36\textwidth}
\captionsetup{size=small}
  \caption{Single Source $(I=1)$ }\label{fig:Effect of mu_a}
  \begin{tikzpicture}[scale=0.70, transform shape]
\begin{axis}[
    xlabel={\large \# of revisions ($\mu$)},
    y tick label style={font=\large},
    x tick label style={font=\large},
    ylabel={\large VAMS (\%) },
    xtick={0,1,2,3,4,5,6,7,8,9,10},
    legend pos=south east,
    ymajorgrids=false,
    grid style=dashed,
]

\addplot[
    color=color1,
    mark=*,
    ]
    coordinates {
    (0,0) (1,80.03668525) (2,97.73911843) (3,100) (4,100)
    };
    
\addplot[
    color=color2,
    mark=*,
    ]
    coordinates {
(0,0) (1,66.92804602) (2,90.21512147) (3,98.56481985) (4,100) (5,100)
    };

\addplot[
    color=color3,
    mark=*,
    ]
    coordinates {
    (0,0) (1,57.50357029) (2,84.51789857) (3,95.09611038) (4,99.21310004) (5,100) (6,100)
    };    
\addplot[
    color=color4,
    mark=*,
    ]
    coordinates {
    (0,0) (1,48.22810898) (2,77.08077653) (3,90.7273416) (4,96.49230093) (5,99.49299934) (6,100) (7,100)
    };

\addplot[
    color=color5,
    mark=*,
    ]
    coordinates {
 (0,0) (1,42.40666731) (2,72.34787595) (3,86.1874389) (4,93.98585874) (5,97.69074738) (6,99.64384894) (7,100) (8,100)
    };

\addplot[
    color=color6,
    mark=*,
    ]
    coordinates {
      (0,0) (1,38.40758516) (2,66.4832187) (3,80.97464332) (4,90.44748525) (5,95.12400533) (6,97.93867263) (7,99.70836967) (8,100) (9,100)
    };

\end{axis}
\end{tikzpicture} 
\end{subfigure}
\begin{subfigure}{.40\textwidth}
\captionsetup{size=small}
   \caption{Multi Source $(I=5)$ }\label{fig:Effect of mu_b}
  \begin{tikzpicture}[scale=0.70, transform shape]
\begin{axis}[
    xlabel={\large \# of revisions ($\mu$)},
    y tick label style={font=\large},
    x tick label style={font=\large},
    xtick={0,1,2,3,4,5,6,7,8,9,10},
    legend pos=south east,
    ymajorgrids=false,
    grid style=dashed,
]
 \addlegendimage{empty legend}
\addplot[
    color=color1,
    mark=*,
    ]
    coordinates {
    (0,0) (1,71.09754335) (2,95.67639368) (3,100) (4,100)
    };

\addplot[
    color=color2,
    mark=*,
    ]
    coordinates {
    (0,0) (1,57.93467544) (2,87.61202046) (3,97.18592295) (4,100) (5,100)
    };    
\addplot[
    color=color3,
    mark=*,
    ]
    coordinates {
   (0,0) (1,50.07557984) (2,77.4843496) (3,91.9385564) (4,98.73765719) (5,100) (6,100)
    };    
\addplot[
    color=color4,
    mark=*,
    ]
    coordinates {
    (0,0) (1,42.54324334) (2,67.02859533) (3,84.40766641) (4,95.19237909) (5,98.89899278) (6,100) (7,100)
    };

\addplot[
    color=color5,
    mark=*,
    ]
    coordinates {
      (0,0) (1,36.81671619) (2,58.78494896) (3,77.29852046) (4,88.35894798) (5,95.9021184) (6,99.01497674) (7,100) (8,100)
    };

\addplot[
    color=color6,
    mark=*,
    ]
    coordinates {
    (0,0) (1,32.36552789) (2,53.75174794) (3,70.51598578) (4,83.1592124) (5,90.80655221) (6,96.06234542) (7,98.84751827) (8,100) (9,100)
    };
    \addlegendentry{\hspace{-.75cm}\textbf{Stages ($T$)}}
    \addlegendentry{5}  
    \addlegendentry{6} 
    \addlegendentry{7}
    \addlegendentry{8}
    \addlegendentry{9} 
    \addlegendentry{10} 
\end{axis}
\end{tikzpicture} 
\end{subfigure}
\end{figure}
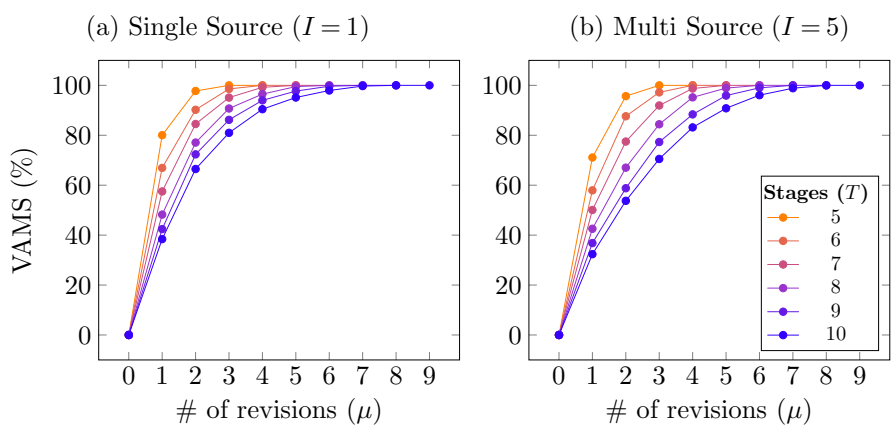

We examine the potential of AMSP in achieving results comparable to those achieved under full flexibility, despite the presence of decision commitment. Figure~\ref{fig:Effect of mu} shows the effect of increasing the number of revisions $\mu$ on the VAMS. Each point represents the VAMS for different values of $\mu$, while the lines represent varying planning horizon lengths $T$. In Figures \ref{fig:Effect of mu_a} and \ref{fig:Effect of mu_b}, we explore two cases: one with a single production source ($I = 1$), and the other with multiple production sources ($I = 5$). In each case, the planning horizon length $T$  varies from 5 to 10, and the number of revisions $\mu$ is within the range from $0$ to $T-1$.

\begin{insight}\label{LS:Insight1}
For a fixed planning horizon $T$, the cost-benefit of increasing the degree of allowed flexibility (i.e., $\mu$) on the Value of the AMSP is diminishing. This implies that in situations where decision commitment is needed in practice, we only need a small number of revision points to obtain most of the value of the fully flexible classic multistage stochastic programming approach. 
\end{insight}

In Figures~\ref{fig:Effect of mu_a} and \ref{fig:Effect of mu_b}, we observe that as the number of revisions increases, the VAMS improves for all values of $T$. However, this improvement diminishes as the number of revisions continues to increase. 
This finding is promising, as it demonstrates that AMSP, even with limited revisions, can achieve comparable results to a fully flexible multistage stochastic programming case.
For instance, we obtain a VAMS exceeding 90\%, if we revise the decisions twice $(\mu = 2)$ in case of a six-stage planning horizon $(T=6)$. This also holds if we revise decisions three times in an eight-stage planning horizon or four times in a ten-stage planning horizon. 
Furthermore, in specific cases, the VAMS reaches $\%100$ without the need to full flexibility. 



\subsubsection{Importance of Multiple Revisions.} 
\begin{figure}[h]
\small
    \caption{Single versus multiple revision points }
    \label{fig: Effect of Stage}
\begin{subfigure}{.34\textwidth}
\captionsetup{size=small}
  \caption{$\mu$=1, \citep{basciftci2019adaptive}}\label{fig: Effect of Stage_a}
  \begin{tikzpicture}[scale=0.68, transform shape]
\begin{axis}[
    xlabel={\large \# Stages ($T$) },
    y tick label style={font=\large},
    x tick label style={font=\large},
    ylabel={\large VAMS (\%) },
    ymin=0, ymax=110,
    xtick={0,20,40,60,80,100},
    xtick={4,5,6,7,8,9,10},
    legend pos=north east,
    legend style={legend title={Sources}},
    ymajorgrids=false,
    grid style=dashed,
]

\addplot[
    color=color1,
    mark=*,
    ]
    coordinates {
(4, 91.51693956)
(5, 80.03668525)
(6, 66.92804602)
(7, 57.50357029)
(8, 48.22810898)
(9, 42.40666731)
(10, 38.40758516)
    };

\addplot[
    color=color3,
    mark=*,
    ]
    coordinates {
(4, 90.54652263)
(5, 74.67590319)
(6, 64.65823372)
(7, 55.62771285)
(8, 46.77881808)
(9, 39.78188908)
(10, 34.56867602)
    };

\addplot[
    color=color6,
    mark=*,
    ]
    coordinates {
(4, 91.40151949)
(5, 71.09754335)
(6, 57.93467544)
(7, 50.07557984)
(8, 42.54324334)
(9, 36.81671619)
(10, 32.36552789)
    };


\end{axis}
\end{tikzpicture} 
\end{subfigure}
\begin{subfigure}{.315\textwidth}
\captionsetup{size=small}
  \caption{$\mu=3$, (AMSP - this paper)}\label{fig: Effect of Stage_b}
  \begin{tikzpicture}[scale=0.68, transform shape]
\begin{axis}[
    xlabel={\large \# Stages ($T$)},
    y tick label style={font=\large},
    x tick label style={font=\large},
    ymin=0, 
    xtick={0,20,40,60,80,100},
    xtick={4,5,6,7,8,9,10},
    legend pos=south east,
    ymajorgrids=false,
    grid style=dashed,
]
\addlegendimage{empty legend}
\addplot[
    color=color1,
    mark=*,
    ]
    coordinates {
  (4, 100)
(5, 100)
(6, 98.56481985)
(7, 95.09611038)
(8, 90.7273416)
(9, 86.1874389)
(10, 80.97464332)
    };

\addplot[
    color=color3,
    mark=*,
    ]
    coordinates {
(4, 100)
(5, 100)
(6, 100.0385362)
(7, 95.66478746)
(8, 89.14043453)
(9, 81.81284742)
(10, 74.64431481)
    };

\addplot[
    color=color6,
    mark=*,
    ]
    coordinates {
(4, 100)
(5, 100)
(6, 97.18592295)
(7, 91.9385564)
(8, 84.40766641)
(9, 77.29852046)
(10, 70.51598578)
    };

\end{axis}
\end{tikzpicture} 
\end{subfigure}
\begin{subfigure}{.30\textwidth}
\captionsetup{size=small}
  \caption{$\mu=5$, (AMSP - this paper)}\label{fig: Effect of Stage_c}
  \begin{tikzpicture}[scale=0.68, transform shape]
\begin{axis}[
    xlabel={\large \# Stages ($T$)},
    y tick label style={font=\large},
    x tick label style={font=\large},
    ymin=0, 
    xtick={0,20,40,60,80,100},
    xtick={4,5,6,7,8,9,10},
    legend pos=south east,
    ymajorgrids=false,
    grid style=dashed,
]
\addlegendimage{empty legend}
\addplot[
    color=color1,
    mark=*,
    ]
    coordinates {
(4, 100)
(5, 100)
(6, 100)
(7, 100)
(8, 99.49299934)
(9, 97.69074738)
(10, 95.12400533)
    };

\addplot[
    color=color3,
    mark=*,
    ]
    coordinates {
(4, 100)
(5, 100)
(6, 100)
(7, 100)
(8, 99.38792922)
(9, 97.24727831)
(10, 94.31579773)
    };

\addplot[
    color=color6,
    mark=*,
    ]
    coordinates {
(4, 100)
(5, 100)
(6, 100)
(7, 100)
(8, 98.89899278)
(9, 95.9021184)
(10, 90.80655221)
    };

    \addlegendentry{\hspace{-.75cm}\textbf{Sources (I)}}
    \addlegendentry{1}  
    \addlegendentry{3} 
    \addlegendentry{5}
\end{axis}
\end{tikzpicture} 
\end{subfigure}
\end{figure}
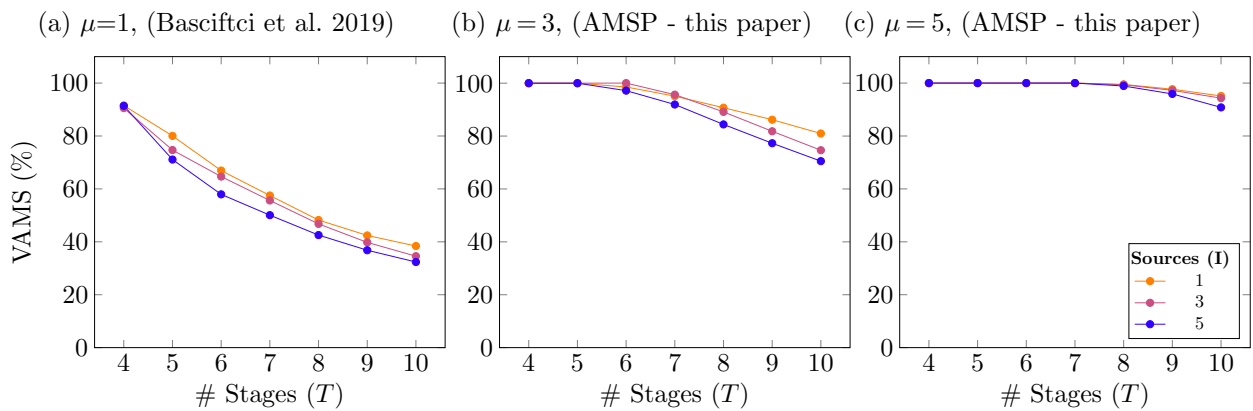

We analyze the potential of considering multiple revision points $(\mu > 1)$ during the planning horizon, compared to considering a single revision point ($\mu = 1$) as introduced by \cite{basciftci2019adaptive}. Figure~\ref{fig: Effect of Stage} show how the VAMS changes as the length of the planning horizon increases under different numbers of sources in the lot-sizing problem. In Figures~\ref{fig: Effect of Stage_a}, \ref{fig: Effect of Stage_b}, and \ref{fig: Effect of Stage_c}, $\mu$ is fixed to 1, 3, and 5, respectively. Each figure includes results for planning horizons with $T$ stages ranging from 4 to 10 and $I$ sources varying among 1, 3, and 5. 

\begin{insight} \label{LS:Insight2}
AMSP offers a clear benefit compared to adaptive two-stage stochastic programming. Whereas adaptive two-stage stochastic programming falls short in the case of relatively long planning horizons, AMSP remains of high quality. 
\end{insight}
 
Zooming in on Figure~\ref{fig: Effect of Stage_a}, we observe a sharp decline in the VAMS as the planning horizon lengths increase. For instance, the spread of the VAMS, across a varying number of sources, ranges between \%90 to \%92 in a four-stage planning horizon, whereas this range is between \%32 and \%38 for a ten-stage planning horizon. Moreover, the VAMS is likely to converge to 0 for even longer planning horizons, turning into a static approach asymptotically. Comparing this against AMSP in case $\mu = 3$, as shown in Figure~\ref{fig: Effect of Stage_b}, we see that the VAMS remains almost 100\% until the planning horizons with six stages. More notable, the spread of VAMS is between \%70 to \%82 for the 10-period planning horizon, significantly higher than the single revision case in Figure~\ref{fig: Effect of Stage_a}. Further increasing the number of revision points makes the VAMS notably higher for longer planning horizons, as seen in  Figure~\ref{fig: Effect of Stage_c}.

Our findings show that VAMS is a reliable approach toward long-term decision-making and highlight the crucial importance of allowing multiple and different revision points for each decision variable in achieving high-quality decisions. In other words, revising decisions only once - as recommended by adaptive two-stage stochastic programming \citep{basciftci2019adaptive} - falls short as the number of stages increases, unless the decision-maker has very limited flexibility by allowing at most one revision during planning. 

\subsubsection{Structure of Optimal Revision Points.}
Optimally choosing revision points is critical in enhancing the efficiency of decision-making. Our observations reveal that the relationship between different revision decisions is not straightforward and requires further analysis. To provide a concrete example, we study an adaptive multistage stochastic uncapacitated lot-sizing problem with a single source and five stages. We enumerate all possible combinations of revision points for different values of $\mu$ and solve Model \eqref{model:MLS} under these revision points. We note that enumeration for instances with multiple sources and longer planning horizons is computationally intractable.

\begin{table}[h]
\small
      \caption{VAMS for single source lot-sizing problem under fixed revisions } \label{tab:LS_Revisions}
      \centering
        \begin{tabular}{cr|cr|cr}
        \toprule
       \multicolumn{2}{c}{$\mu=1$} & \multicolumn{2}{|c|}{$\mu=2$} &\multicolumn{2}{c}{$\mu=3$} \\   \midrule
                 Revision points & VAMS (\%) &  Revision points & VAMS (\%)   & Revision points & VAMS (\%)\\  \midrule                 
        2 &  7  & 2, 3 & 78   & 2, 3, 4 & 79 \\ 
        3 &  78 & 2, 4 & 79   & \textbf{2, 3, 5} & \textbf{100}\\ 
        \textbf{4} &  \textbf{79} & 2, 5 & 59   & \textbf{2, 4, 5} & \textbf{100}\\ 
        5 &  59 & 3, 4 & 79   & \textbf{3, 4, 5} & \textbf{100}\\ 
          &     & \textbf{3, 5} &  \textbf{100}  &       & \\ 
          &     & \textbf{4, 5}& \textbf{100}  &       &  \\ \bottomrule
\end{tabular}
\end{table}

\begin{insight}
    Selecting the optimal revision stages is essential. 
    Revising decisions at the optimal stages can outperform revising decisions more often in arbitrary stages.
\end{insight}

Table \ref{tab:LS_Revisions} reports the VAMS for each possible combination of revision points, and indicates the best ones in bold. For $\mu=1$, we observe that fixing the revision point to the 4th stage yields the highest VAMS value while the VAMS in the 3rd stage is only marginally lower than the optimal value. However, the VAMS significantly decreases when revisions are made in the 2nd or 5th stages. When two revisions are allowed, the best results are achieved with revision points (3,5) and (4,5). Since the 4th and 3rd periods are the most effective ones in the one revision case, the optimal combination for the two-revision case can be anticipated as (3,4). However, our analysis indicates that this is not the case, highlighting the complex nature of revision decisions. 
The results show the importance of not only the number but also the timing of revisions. For instance, revising twice in (3,5) or (4,5) is more effective than three times in (2,3,4). Similarly, revising once in 3 or 4 is more effective than twice in (2,5). Overall, the result of the adaptive method is highly impacted by the timing of revision, emphasizing the significance of optimally choosing the revision points in situations with limited flexibility.

\subsection{Adaptive Multistage Stochastic Generation Expansion Planning Problem}\label{Sec:GenerationExpansion}
We study a generation expansion planning (GEP) problem: a process used by energy system planners to determine the optimal mix of generation sources to meet future power demand. This problem is commonly tackled using multistage stochastic programming methods, where the goal is to primarily adjust the generation expansion decisions (i.e., the installed capacity of energy generation sources) and consequently the operational decisions (i.e., how much power to produce) as uncertainties realize over a multi-period planning horizon. However, fully flexible policies obtained by multistage stochastic programming are not practically feasible for expansion decisions due to restrictions such as contractual commitments, and lead time requirements \citep{munoz2013engineering,careri2011generation}. In such cases where flexibility is limited, an effective solution is selecting critical stages for revising decisions while committing to these decisions in between. Thus, we introduce an AMSP approach to solve the GEP problem. We consider the allowed -- practically deemed feasible -- flexibility level of the underlying GEP process. Consequently, expansion decisions will only be revised at the most critical stages determined by our optimization model, leading to more practically efficient planning outcomes.

We consider a scenario tree $\mathcal{T}$ with a finite planning horizon $\mathcal{P}=\{1, \dots, T\}$. The probability of uncertainty realization at node $n \in \mathcal{T}$ is denoted by $p_n$. Each stage $t$ is divided into subperiods $\mathcal{K}_{t}$ reflecting hours within that stage. Then, the total number of hours at node $n \in \mathcal{T}$ in subperiod $k \in \mathcal{K}_{t_n}$ is denoted as $h_{nk}$. We consider a set of different generators denoted by $i \in \mathcal{G}$. The initial number of generators of type $i \in \mathcal{I}$ at the beginning of the planning is denoted as $x^0_{i}$. The construction of a generator of type $i \in \mathcal{I}$ at node $n \in \mathcal{T}$ incurs a capital cost of $c^c_{in}$, a fixed operations and maintenance (O\&M) cost of $c^f_{in}$, and a variable O\&M cost of $c^v_{in}$. The maximum number of generators of type $i \in \mathcal{I}$ that can be constructed at node $n \in \mathcal{T}$ is limited to $x_{i}^{\max}$. Each generator type $i \in \mathcal{I}$ has a maximum generation capacity of $m_i^{\max}$. Additionally, each generator $i$ in $\mathcal{I}$ has a capacity factor $l_{ink}$, which indicates the amount of energy it can produce at node $n \in \mathcal{T}$ during subperiod $k \in \mathcal{K}_{t_n}$, relative to its maximum generation capacity. The fuel price of a generator of type $i \in \mathcal{I}$ at node $n \in \mathcal{T}$ is denoted as $c^u_{in}$. The hourly demand at node $n \in \mathcal{T}$ in subperiod $k \in \mathcal{K}_{t_n}$ is given by $d_{nk}$. Failure to meet the demand incurs a penalty cost of $c^p$. The yearly interest rate is represented by $\gamma$. Note that capacity factor  $l_{ink}$ and demand $d_{nk}$ are stochastic parameters. 

Let $x_{in}$ be an integer decision variable that indicates the number of type $i \in \mathcal{I}$ generators constructed at node $n \in \mathcal{T}$. Let $y_{i n k}$  be a continuous decision variable for the generation of a type $i \in \mathcal{I}$ generator in subperiod $k\in \mathcal{K}_{t_{n}}$ at node $n \in \mathcal{T}$. Finally, let $u_{n k}$  be a continuous variable for unsupplied demand in subperiod $k\in \mathcal{K}_{t_{n}}$ at node $n \in \mathcal{T}$. We then formulate the adaptive multistage GEP based on the formulation in \cite{zou2018partially}, as follows:
\begin{subequations}\label{model:GE}
\begin{align}
\min \quad & \sum_{n \in \mathcal{T}} p_{n} \frac{1}{(1 + \gamma)^{t_{n} - 1}}  \Biggl(\sum_{i \in \mathcal{G}} \ \left(c^c_{in} + \sum_{t=t_n}^{T} \frac{c^f_{in}}{(1 + \gamma)^{t - t_{n}}}  \right) m_i^{\max} x_{i n}  \notag \\
& 
\hspace{3.3cm}  +\sum_{i \in \mathcal{G}} \sum_{k \in \mathcal{K}_{t_{n}}}  (c^v_{in} + c^u_{in}) h_{nk} y_{ink} + \sum_{k \in \mathcal{K}_{t_{n}}}  c^p h_{nk}  u_{ n k} \Biggr)
\label{GE_obj} \\
\text{s.t.} \quad &  \sum_{m \in P(n)} x_{i m} + x^0_{i} \geq  y_{i n k} \hspace{4.05cm}  \forall i \in \mathcal{G}, k \in \mathcal{K}_{t_{n}}, n \in \mathcal{T} 
 \label{GE_gen}  \\
&   x_{i n}  \leq x_{in}^{\max } \hspace{6cm} \forall i \in \mathcal{G}, n \in \mathcal{T}\label{GE_cap}  \\
&  \sum_{i \in \mathcal{G}} l_i y_{i n k} + u_{n k} \geq d_{n k} \hspace{4.1cm} \forall k \in \mathcal{K}_{t_{n}}, n \in \mathcal{T}   \label{GE_demand}\\
&  x_{i n} \in \mathbb{Z}_{+}, y_{i n k} \in \mathbb{R}_{+} \hspace{4.55cm} \forall i \in \mathcal{G}, k \in \mathcal{K}_{t_{n}}, n \in \mathcal{T} \label{GE_set} \\
& \eqref{AdaptiveModel_f}-\eqref{AdaptiveModel_g}, \eqref{combinedNACs} \notag
\end{align}\end{subequations}
The objective function \eqref{GE_obj} minimizes the expected cost of generator construction, fixed and variable operations and maintenance, fuel prices, and penalty for unsupplied demand over a planning period, discounting it to the beginning of the planning period. Constraint  \eqref{GE_gen} guarantees that the generation amount does not exceed the total generation capacity of constructed generators. Constraint \eqref{GE_cap} limits the number of generators that can be constructed. Constraint \eqref{GE_demand} ensures demand is satisfied. Finally, constraints \eqref{AdaptiveModel_f} - \eqref{AdaptiveModel_g}, \eqref{combinedNACs} are the necessary NACs.

\subsubsection{Experimental Setup.}
We set model parameters based on real-world data and optimize the yearly investment choices for various energy sources, ensuring each source operates within its capacity limits while meeting the system's demand for each sub-period. We consider 2020 as the base year. We include five different types of generators: a traditional gas-powered combined-cycle generator, a gas-powered combined-cycle generator with carbon capture (CC), and three renewable energy sources: onshore wind, offshore wind, and solar PV. We set the capital and O\&M cost parameters based on the ``Capital Cost and Performance Characteristic Estimates for Utility Scale Electric Power Generating Technologies" report from the Energy Information Administration (EIA) \citep{Capital2020}. The report provides various characteristics of these technologies, such as their net nominal capacity, capital cost, and fixed and variable O\&M costs. We set the cost parameters at the base year ($c^c_{i0}$, $c^f_{i0}$, $c^v_{i0}$) according to values provided in Table \ref{Table:GEP Data}.  By utilizing the ``EIA Electric Power Annual 2021" \citep{Electric2022}, we set the natural gas price for the combined-cycle power plants, denoted as $c^u_{in}$, to 0.04\$  per kWh.
Based on \cite{basciftci2019adaptive}, we assume a yearly decrease of 10\% in capital costs for renewable energy sources and 5\% for the combined-cycle with CC, reflecting expected technology development. Additionally, we consider a 10\% yearly increase in variable O\&M costs and fuel prices for both the combined-cycle and the combined-cycle with CC. Other cost parameters remain constant throughout the planning horizon.

\begin{table}[h]
\small
\caption{Generator costs in the baseline year 2020}\label{Table:GEP Data}
\centering
\begin{tabular}{l|rrrr }
\toprule
        Type of Generator & Nominal Capacity & Capital &	Fixed O\&M 
& Variable O\&M  \\
         & (kW) &  (\$/kW) &	(\$/kW-year)
&  (\$/MWh) \\

        \midrule
       Combined-Cycle & 418 & 1084.0	&14.1	&2.6   \\    
       Combined-Cycle with 90\% CC &377 & 2481.0 &	27.6	&5.8   \\    
       Onshore Wind &200 & 1265.0 &	26.4&	0 \\    
       Offshore Wind &400 & 4375.0 &	110.0 &	0  \\  
       Solar PV &150 & 1313.0 & 15.3 &	0   \\ \bottomrule  
\end{tabular}
\end{table}



Following the study by \cite{min2018long}, we assume that each year has four subperiods differentiated by their demand level: base, off-peak, shoulder, and peak. The demand per subperiod ($d_{nk}$) is calculated by multiplying the node's demand ($d_n$) by the corresponding weights ($w_k$) of 0.9, 1.1, 1.3, and 1.5, respectively, representing the base, off-peak, shoulder, and peak subperiods. The subperiods have different durations within a year, with the base period covering 55\%, the off-peak period covering 40\%, the shoulder period covering 4.95\%, and the peak period covering 0.05\%. Based on \cite{singh2009dantzig}, the scenario tree is generated to capture the stochastic nature of the demand and capacity factor for renewable energy sources as in Algorithm \ref{ScenarioTreeGeneration} in Appendix \ref{sec:AppendixScenarioTreeAlg}.

We assume that demand at the root node $d_1$ is 1000MwH. Each node $n$ of the scenario tree, except the root node, has a randomly allocated growth factor $\lambda^d_{nk}$ for all subperiods $k \in \mathcal{K}_{t_n}$. These growth rates are derived from a normal distribution with a mean of 5\% ($\mu^d$) and a standard deviation of 5\% ($\sigma^d$), reflecting the yearly increase in global electricity demand over recent years as reported in \cite{IEA2023}. The amount of energy produced by renewable energy sources is determined by their capacity factor $l_{ink}$. To reflect the uncertainty of renewable energy production, we generate random samples of capacity factors from a normal distribution with $\mathcal{N}(30\%,10\%)$ for onshore wind, $\mathcal{N}(60\%,5\%)$ for offshore wind, and $\mathcal{N}(20\%,10\%)$ for solar PV. 
We assume the network consists of only combined-cycle generators in the base year. Accordingly, we adjust the initial capacities $x_i^0$ to ensure that they can cover the average demand if operated at full capacity. We set the yearly interest rate $\gamma$ to 5\%.



\subsubsection{Optimal Revision Points.}

Table \ref{Table:GEP_Revisions} shows the optimal decisions for the revision points and VAMS for each type of generator across different values of $\mu$ ranging from 1 to 4 when $T = 9$.

\begin{table}[h]
\small
      \centering
           \caption{Optimal revision points for the GEP with $T=9$}\label{Table:GEP_Revisions}
          
        \begin{tabular}{l|rrrr  } \toprule
         {Type of Generator} & \multicolumn{4}{c}{Number of revisions} \\ 
         & $\mu=$1 & $\mu=$2 & $\mu=$3  & $\mu=$4 \\ \midrule 
       Combined-Cycle & 4 & 4, 7 & 2, 4 ,6 & 2, 4, 6, 7 	  \\ 
       Combined-Cycle with 90\% CC & 9 & 6, 8 & 5, 6, 9  & 5, 7, 8, 9 \\   
       Onshore Wind & 5 & 4, 7 & 4, 6, 8  & 4, 6, 8, 9 \\   
       Offshore Wind & 7 & 6, 8 &	6, 8, 9  & 4, 6, 8, 9 \\ 
       Solar PV & 9 & 7, 9 & 7, 8, 9 &  5, 7, 8, 9  \\ \midrule
       \textbf{VAMS (\%)} &  76.7 & 88.3 & 92.9   & 95.7 \\ \bottomrule 
\end{tabular}
\end{table}

\begin{insight}
Optimal revision points for the traditional combined-cycle gas-powered generator are relatively early compared to the revision points for renewable energy sources and the combined-cycle gas-powered generator with CC. This shows the importance of allowing each generator to be revised at different stages, as advocated by AMSP.
\end{insight}

The structure of optimal revision points can be explained by the underlying dynamics of the GEP. Capacity expansion decisions for combined-cycle generators occur in earlier periods due to increasing O\&M costs and natural gas prices over the years. Consequently,  optimal revision points for combined-cycle generators are at earlier periods. However, as capital costs of the combined-cycle with 90\% CC decrease over time due to advancements in CC technologies, expansions are shifted to later periods. Likewise, renewable capacities also expand during these later periods. The main reason is the declining capital costs of renewables and the increased variability of their capacity factors within the scenario tree. As a result, capacity expansion amounts during these later periods show notable differences depending on the realization of uncertainties. Therefore, optimal revision points for wind and solar are at later periods. To summarize, optimizing power generation expansion planning requires determining tailored revision points for each generator type, considering evolving costs, technological advancements, and market dynamics.

The results indicate a positive correlation between the number of revisions and VAMS. Moreover, as the number of revisions increases, we observe a diminishing return, which aligns with Insight \ref{LS:Insight1}.  Additionally, our findings support Insight \ref{LS:Insight2}, as multiple revisions yield significantly higher VAMS values than a single revision. 




\subsubsection{Computational Performance.}
We analyze the impact of our NACs elimination techniques and our solution algorithm on the computational performance of AMSP. To this end, we first replace the NACs in Model \eqref{model:GE} (corresponding to formulation  $(AMS_\mu+$)) with NACs in formulation $(AMS_{\mu})$. Accordingly, we obtain model $\{\eqref{GE_obj} : \eqref{GE_gen}-\eqref{GE_set}, \eqref{AdaptiveModel_d}-\eqref{AdaptiveModel_g} \}$, and we call it \texttt{GEP} which corresponds to the original AMSP formulation without the NACs elimination techniques. We call the improved model, Model \eqref{model:GE}, as \texttt{GEP+}. 

\begin{table}[h]
\small
    \centering
     \caption{Performance of \texttt{GEP}, \texttt{GEP+} and Algorithm \ref{SolutionAlgorithm}}\label{tab:GEP_Computational3}
    \begin{tabular}{ll|rrrrrr}
    \toprule
         &  & \multicolumn{2}{c}{\texttt{GEP}} & \multicolumn{2}{c}{\texttt{GEP+} } & \multicolumn{2}{c}{Algorithm \ref{SolutionAlgorithm}  } \\ \cmidrule(lr){3-4} \cmidrule(lr){5-6} \cmidrule(lr){7-8}
        Stage & $\mu$ & Time (s) & VAMS (\%) & Time (s) & VAMS (\%) & Time (s)  & VAMS (\%) \\ \midrule
         5 & 1 & 5.7  & 75.2 & 4.1  & 75.2 & 0.7  & 75.2 \\
         & 2 & 1.7  & 88.2 & 0.8  & 88.2 & 0.5  & 88.2 \\ 
         & 3 & 0.6  & 97.7 & 0.5  & 97.7 & 0.2  & 97.7 \\ 
        6 & 1 & 188.6  & 78.0 & 6.1  & 78.0 & 3.6   & 78.0 \\ 
         & 2 & 16.4  & 90.1 & 8.2  & 90.1 & 4.3  & 90.1 \\
          & 3 & 10.6  & 95.3 & 6.6  & 95.3 & 2.8  & 95.3 \\ 
          & 4 & 7.8 & 98.8 & 2.4  & 98.8 & 2.4  & 98.8 \\ 
        7 & 1 & -  & - & 70.6  & 78.4 & 21.4  & 78.4 \\
        ~ & 2 & 389.5  & 90.7 & 67.0  & 90.7 & 14.6  & 90.7 \\ 
        ~ & 3 & 113.7  & 94.4 & 26.7  & 94.4 & 11.5  & 94.4 \\ 
         ~ & 4 & 92.3  & 98.2 & 18.7  & 98.2 & 8.4  & 98.2 \\
        8 & 1 & -  & - & 692.7  & 77.9 & 203.1  & 77.9 \\ 
        ~ & 2 & -  & - & 264.5  & 88.9 & 200.0  & 88.9 \\
        ~ & 3 & -  & - & 136.2  & 94.3 & 86.1  & 94.3 \\ 
        ~ & 4 & -  & - & 88.8  & 96.7 & 50.5  & 96.7 \\ 
        9 & 1 & -  & - & 18000.0  & 76.5 & 1061.6  & 76.7 \\ 
        ~ & 2 & -  & - & 12605.7  & 88.3 & 2573.5  & 88.3 \\ 
        ~ & 3 & -  & - & 2552.7  & 92.9 & 1596.8  & 92.9 \\ 
         ~ & 4 & -  & - & 1395.6  & 95.7 & 433.7  & 95.7 \\ \midrule
        Average &  &  -  & - & 1892.0   & 89.3 & 331.7   & 89.3 \\ \bottomrule
    \end{tabular}
\end{table}

We present the computational results of \texttt{GEP}, \texttt{GEP+}, and Algorithm \ref{SolutionAlgorithm} in Table \ref{tab:GEP_Computational3}, which includes solution times, optimality gaps, and VAMS for each approach under various instances. These results are under a threshold of 0.1\% optimality gap and a time limit of 5 hours. For our analysis, we generate instances with varying flexibility levels ($\mu$ ranging from 1 to 4) and numbers of stages ($T$ ranging from 5 to 9). The scenario tree utilized in our experiments consists of 3 branches per node. 

We observe that \texttt{GEP} is outperformed by \texttt{GEP+} in all instances. \texttt{GEP} gives memory error failing to solve any instance with more than eight stages and does not provide feasible solutions within the time limit. Moreover, the instance with seven stages and one revision cannot be solved. On the other hand, \texttt{GEP+} can solve all but one instance to optimality. The optimality gap of that single instance (with nine stages and one revision) is 0.5\%. The results show that including the elimination techniques of the NACs leads to an average decrease in computation time of 59\%, with up to 97\%. We also observe that as the size of the problem increases, the superiority of \texttt{GEP+} over \texttt{GEP} in terms of solution time becomes more prominent. This is due to the significant increase in the number NACs for larger instances as reported in Table \ref{tab:Number of NACs}.

We show that the solution time can be further reduced with Algorithm \ref{SolutionAlgorithm}. Compared to \texttt{GEP}, the solution time is decreased on average by 82\% with up to 98\%. Compared to \texttt{GEP+}, solution time is decreased on average 56\% with up to 94\%. It is worth noting that although Algorithm \ref{SolutionAlgorithm} includes heuristic steps, it can still attain the optimal value for all instances. Overall, Algorithm \ref{SolutionAlgorithm} proves to be highly effective in reducing solution times and delivering optimal solutions for AMSPs.

\section{Conclusions}\label{sec:Conclusions}
In this paper,  we propose a new optimization approach called AMSP to achieve the optimal  balance between decision flexibility and commitment by deciding on the best times to revise decisions under uncertainty. We present a generic mathematical formulation for the proposed approach and prove the NP-Hardness of the proposed problem. Then, we derive theoretical properties for elimination of NACs as their existence increases the problem complexity. We provide a tailored solution algorithm for solving the generic class of AMSPs with mixed-integer master problem and subproblems, 
by combining L-shaped and Benders decomposition methods. To further accelerate the algorithm, we include a heuristic-based preprocessing step that can eliminate a set of candidate feasible solutions.



Computational experiments on lot-sizing and generation expansion planning  demonstrate the significant benefits of optimally choosing revision moments when flexibility is limited. By making revisions only at these optimal moments, businesses can produce results comparable to a fully flexible setting while respecting the flexibility limitations of the decision-making process. Remarkably, in some cases, we show that the optimal timing of revisions can even outperform frequency. Next to these practical insights, we demonstrate that eliminating NACs improves computational performance over various instances by an average of 59\%, with our solution algorithm providing a further 56\% enhancement on average.

\section*{Acknowledgements}
This project has received funding from the Fuel Cells and Hydrogen 2 Joint Undertaking (now Clean Hydrogen Partnership) under Grant Agreement No 875090. This Joint Undertaking receives support from the European Union's Horizon 2020 research and innovation programme, Hydrogen Europe and Hydrogen Europe Research. Albert H. Schrotenboer is supported by a VENI research talent grant from the Dutch Science Foundation (NWO).

\bibliographystyle{apalike}
\bibliography{references}

\begin{APPENDICES} 

\section{Proofs of Lemmas and Propositions}
This section presents the proofs of the lemmas and propositions discussed in the main paper. To enhance readability, corresponding statements are restated below. 

\subsection{Proof of Lemma \ref{Lem:binding}} \label{App:Binding}

\textbf{Lemma 1} \textit{For $(AMS_\mu)$, there exist an optimal solution where $r_{iT} \leq \mu$ is binding for all $i \in \mathcal{I}$.}

\proof{Proof:}
Assume that constraint $r_{iT}^* \leq \mu$ has a slack $\epsilon \in \{1,\dots,\mu\}$ for at least one $i \in \mathcal{I}$ at an optimal solution $\mathbf{r^*}$. 
Let $\mathcal{I}' \subseteq \mathcal{I}$ be the set of indices of the state variables with slack at this solution. 
Since any such solution under $\mathbf{r^*}$ has an objective function value which is same as the solutions we can construct by imposing $r_{iT} = \mu$ for all $i \in \mathcal{I}'$, the desired result follows as we can obtain a different but still optimal solution where $r_{iT} = \mu$ for all $i \in \mathcal{I}$. 
\Halmos

\endproof

\subsection{Proof of Proposition \ref{prop:TS-MS}} \label{App:TS-MS}

\textbf{Proposition 1} \textit{For $\mu = 0$, formulation $(AMS_\mu)$ is equivalent to two-stage stochastic programming formulation $(2SP)$. For $\mu = T - 1$, formulation $(AMS_\mu)$ is equivalent to multistage stochastic programming formulation $(MSP)$.}

\proof{Proof:} If $\mu = 0$, then $r_{it} = 0$ for all $i = 1, \dots, I, t = 1, \dots, T$. Then \eqref{AdaptiveModel_d} and \eqref{AdaptiveModel_e} reduce to $x_{im} = x_{in}$, which is analogous to constraint \eqref{Two-stageC} of $(2SP)$. On the other hand, for $\mu = T-1$, there exist an optimal solution where $r_{i(t+1)} = r_{it} + 1 $ with $r_{iT}= T-1$  for $i = 1, \dots, I, t = 1, \dots, T-1 $  based on Lemma \ref{Lem:binding}.  As a result, constraints \eqref{AdaptiveModel_d} and \eqref{AdaptiveModel_e} become redundant. Thus, it can be concluded that $(AMS_{(T-1)})$ is equivalent to $(MSP)$. \Halmos 
\endproof 

\subsection{Proof of Proposition \ref{prop:ObjectiveComparison}} \label{App:ObjectiveComparison}

\textbf{Proposition 2} \textit{$z(2SP) =  z(AMS_0) \geq z(AMS_1) \geq \dots \geq z(AMS_{(T-2)}) \geq z(AMS_{(T-1)}) = z(MSP)$. 
}

\proof{Proof:} To demonstrate this result, we show that the feasible region of $(AMS_{(\mu-1)})$ is contained in the feasible region of $(AMS_{\mu})$ for any $\mu \in \{1,\dots, T-1\}$. Note the constraint sets of the problems $(AMS_{(\mu-1)})$ and $(AMS_{\mu})$ are the same except constraint \eqref{AdaptiveModel_f}, where the upper bound on the $ r_{it}$ variables are affected by the change in $\mu$ value for every $i \in \mathcal{I}$ and $t \in \mathcal{P}$. While $AMS_{(\mu-1)}$ has $ r_{iT} \leq \mu-1$, $AMS_{\mu}$ has $ r_{iT} \leq \mu$ where  $r_{iT} \in \mathbb{Z}_{+}$. Since $\mu-1 \leq \mu $, $ z(AMS_{(\mu-1)}) \geq  z(AMS_{\mu})$ for any $\mu \in \{1,\dots, T-1\}$. \Halmos
\endproof

\subsection{Proof of Proposition \ref{prop:NPHard}} \label{App:NPHard}
\textbf{Proposition 3} \textit{
AMSP 
is NP-Hard.}

\proof{Proof:} 
To prove the Proposition \ref{prop:NPHard}, we demonstrate that the feasibility problem associated with the adaptive multistage problem \eqref{AMSLP} is $\mathcal{N} \mathcal{P}$-Complete. To this end, we construct a feasibility problem, referred as P, whose feasible region is defined in \eqref{NPHardModel}.
Then, we show that the subset sum problem, which is known to be $\mathcal{N} \mathcal{P}$-Complete \citep{garey1980computers},
can be reduced to $\mathrm{P}$ in polynomial time. We specify the subset sum problem as follows: Given the non-negative integers 
 $\{w_{il}^{'}\}_{i \in \mathcal{I}^{'}, l\in \mathcal{L}}$ and $W^{'}$, does there exist $\mathcal{H}_l \subseteq \mathcal{I}^{'}$ for all  $l\in \mathcal{L}$ such that $\sum_{l\in \mathcal{L}}\sum_{i \in \mathcal{H}_l} w_{il}^{'}=W^{'}$?
\begin{subequations}\label{NPHardModel}
\allowdisplaybreaks\begin{align}
\min _{\mathbf{\alpha}, \mathbf{\beta}, \mathbf{t^*}} \quad & 0 \label{NPHardModel_obj} \\
\text{s.t.} \quad & \alpha_{i 2}- \beta_{i 2}=0 & i \in \mathcal{I}  \label{NP:b}\\
& \alpha_{i 3}-\beta_{i 3}=\alpha_{i 1} & i \in \mathcal{I} \\
& \alpha_{i 4}- \beta_{i 4}=0 & i \in \mathcal{I} \\
&\alpha_{i 5}-\beta_{i 5}=\alpha_{i 2} & i \in \mathcal{I} \\
& \alpha_{i 6}-\beta_{i 6}=0 & i \in \mathcal{I} \\
& \alpha_{i 7}-\beta_{i 7}=\alpha_{i 2} & i \in \mathcal{I} \\
& \alpha_{i 8}-\beta_{i 8}=0 & i \in \mathcal{I} \\
&\alpha_{i 9}-\beta_{i 9}=\alpha_{i 4} & i \in \mathcal{I} \\
&\alpha_{i 10}-\beta_{i 10}=0 & i \in \mathcal{I} \\
& \alpha_{i 11}-\beta_{i 11}=\alpha_{i 4} & i \in \mathcal{I} \\
&\alpha_{i 12}-\beta_{i 12}=0 & i \in \mathcal{I} \\
&\alpha_{i 13}-\beta_{i 13}=\alpha_{i 4} & i \in \mathcal{I} \\
& \alpha_{i 14}-\beta_{i 14}=0 & i \in \mathcal{I} \\
& \alpha_{i 15}-\beta_{i 15}=\alpha_{i 4} & i \in \mathcal{I} \label{NP:o} \\
&\sum_{i \in \mathcal{I}} (w_{i1} \alpha_{i 1} + w_{i2} \alpha_{i 2} + w_{i4} \alpha_{i 4})  =W & \label{NP:W} \\
&\alpha_{i m}=\alpha_{i n}, \beta_{i m}=\beta_{i n} \quad \forall m, n \in S_t, \quad t<t^*_{i1} & i \in \mathcal{I} \label{NP:NA1} \\
&\alpha_{i m}=\alpha_{i n}, \beta_{i m}=\beta_{i n}  \quad \forall m, n \in S_t \cap \mathcal{T}(l,t^*_{ik+1}), \ l \in S_{ t^*_{ik}}, \  t^*_{ik+1} >  t \geq t^*_{ik} &  i \in \mathcal{I}   \\
&\alpha_{i m}=\alpha_{i n}, \beta_{i m}=\beta_{i n}  \quad \forall m, n \in S_t \cap \mathcal{T}(l), \quad l \in S_{t^*_{i\mu}}, \quad  t \geq t^*_{i\mu} & i \in \mathcal{I}  \\ 
& t^*_{ik} \in \Gamma_i^*  \subseteq  \{2, \dots, T\} & i \in \mathcal{I}, k \in \{1,\dots,\mu\}  \label{NP:NA4}
\end{align}
\end{subequations}

First observe that the problem \eqref{NPHardModel} is an instance of the adaptive multistage problem by letting 
$ \mathcal{J}=\emptyset$, and $x_{i n}=\left[\begin{array}{c}\alpha_{i n} \\ \beta_{i n}\end{array}\right]$. We consider the scenario tree $\mathcal{T}$ with $T=4$ stages and $B=2$ branches. Each stage of the scenario tree consists of the following nodes: $S_1 = \{1\}, \> S_2 = \{2,3\}, \>  S_3 = \{4,5,6,7\}, \> S_4 = \{8,9,10,11,12,13,14,15\}$. We let $\mathcal{I}=\mathcal{I}^{'}, W=W^{'}$, and $w_{il}=w_{il}^{'}$ for all $i \in \mathcal{I}, l \in \mathcal{L}=\{1,2,4\}$. We note that constraints \eqref{NP:b}--\eqref{NP:W} correspond to constraint \eqref{eq:MScons}. 
Similarly, constraints \eqref{NP:NA1}--\eqref{NP:NA4} refer to the constraints \eqref{AMSLP_d}--\eqref{AMSLP_g}, respectively. In Lemma \ref{Lem: NP-Hard}, we list all possible combinations of revision points for varying $\mu$ from zero to three. We note that the first case, where $\Gamma_i^*=\{\}$, corresponds to a scenario with no revisions during planning. From the second to the fourth case, there is one revision, from the fifth to the seventh case, there are two revisions, and the last case is fully flexible.
\begin{lem}\label{Lem: NP-Hard}
Under the given set of revision points $\Gamma_i^*$ for all $i \in \mathcal{I}$, we can construct the following feasible solutions for the problem \eqref{NPHardModel}: \\
1. If $\Gamma_i^*=\{\}$,  then $\alpha_{i 1}=\alpha_{i 2}=\alpha_{i 4}=0$ for all $i \in \mathcal{I}$.\\
2. If $\Gamma_i^*=\{2\}$,  then $\alpha_{i 1}=1, \alpha_{i 2}=\alpha_{i 4}=0$ for all $i \in \mathcal{I}$.\\
3. If $\Gamma_i^*=\{3\}$, then $\alpha_{i 2}=1, \alpha_{i 1}=\alpha_{i 4}=0 $ for all $i \in \mathcal{I}$.\\
4. If $\Gamma_i^*=\{4\}$, then $\alpha_{i 4}=1,  \alpha_{i 1}=\alpha_{i 2}=0$ for all $i \in \mathcal{I}$.\\ 
5. If $\Gamma_i^*=\{2,3\}$, then $\alpha_{i 1}=\alpha_{i 2}=1, \alpha_{i 4}=0$ for all $i \in \mathcal{I}$.\\
6. If $\Gamma_i^*=\{2,4\}$,  then $\alpha_{i 1}=\alpha_{i 4}=1, \alpha_{i 2}=0  $ for all $i \in \mathcal{I}$.\\
7. If $\Gamma_i^*=\{3,4\}$, then $\alpha_{i 2}=\alpha_{i 4}=1, \alpha_{i 1}=0 $ for all $i \in \mathcal{I}$.\\
8. If $\Gamma_i^*=\{2,3,4\}$, then $\alpha_{i 1}=\alpha_{i 2}=\alpha_{i 4}=1 $ for all $i \in \mathcal{I}$.
\end{lem}

Given a solution $\{\mathcal{H}_l\}_{l \in \mathcal{L}}$ to the subset-sum problem, we can construct a feasible solution for problem \eqref{NPHardModel} based on Lemma \ref{Lem: NP-Hard}. For all $l \in \{1,2,4\}$, we set $\alpha_{i l}=1$ if $i \in \mathcal{H}_l$ and 0 otherwise. For example, in the first case: if $\alpha_{i1}=\alpha_{i2}=\alpha_{i4}=0$, then  $\alpha_{in}=\beta_{in}=0$ for all $i \in \mathcal{I}, n \in \mathcal{T}$. The second case: $\alpha_{i1}=1, \alpha_{i2}=\alpha_{i4}=0$, then $ \alpha_{i3} = 1$, $\alpha_{in}=0$ for all $i \in \mathcal{I}, n \in \mathcal{T}\backslash\{1,3\}$ and $\beta_{in}=0$ for all $i \in \mathcal{I}, n \in \mathcal{T}$ where $t^*_{i1}=2$. The third case: $\alpha_{i2}=1, \alpha_{i1}=\alpha_{i4}=0$, then $ \alpha_{i3}=\alpha_{i5}=\alpha_{i7}=1, \alpha_{in}=0$ for all $i \in \mathcal{I}, n \in \mathcal{T}\backslash\{2,3,5,7\}$ and $ \beta_{i2}=\beta_{i3}=1, \beta_{in}=0$ for all $i \in \mathcal{I}, n \in \mathcal{T}\backslash\{2,3\}$ where $t^*_{i1}=3$. The fourth case: $\alpha_{i4}=1, \alpha_{i1}=\alpha_{i2}=0$, then $ \alpha_{i5}=\alpha_{i6}=\alpha_{i7}=\alpha_{i9}=\alpha_{i11}=\alpha_{i13}=\alpha_{i15}=1, \alpha_{in}=0$ for all $i \in \mathcal{I}, n \in \mathcal{T}\backslash\{4,5,6,7,9,11,13,15\}$ and $\beta_{i4}=\beta_{i5}=\beta_{i6}=\beta_{i7}=1, \beta_{in}=0$ for all $i \in \mathcal{I}, n \in \mathcal{T}\backslash\{4,5,6,7\}$ where $t^*_{i1}=4$. The fifth case: $\alpha_{i1}=\alpha_{i2}=1, \alpha_{i4}=0$, then $ \alpha_{i3}=2,\alpha_{i5}=\alpha_{i7}=1, \alpha_{in}=0$ for all $i \in \mathcal{I}, n \in \mathcal{T}\backslash\{1,2,3,5,7\}$ and $\beta_{i2}=\beta_{i3}=1, \beta_{in}=0$ for all $i \in \mathcal{I}, n \in \mathcal{T}\backslash\{2,3\}$ where $t^*_{i1}=2, t^*_{i2}=3$. The sixth case: $\alpha_{i1}=\alpha_{i4}=1, \alpha_{i2}=0$, then $ \alpha_{i3}=\alpha_{i5}=\alpha_{i6}=\alpha_{i7}=\alpha_{i9}=\alpha_{i11}=\alpha_{i13}=\alpha_{i15}=1, \alpha_{in}=0$ for all $i \in \mathcal{I}, n \in \mathcal{T}\backslash\{1,3,4,5,6,7,9,11,13,15\}$ and $\beta_{i4}=\beta_{i5}=\beta_{i6}=\beta_{in7}=1, \beta_{in}=0$ for all $i \in \mathcal{I}, n \in \mathcal{T}$  where $t^*_{i1}=2, t^*_{i2}=4$. The seventh case: $\alpha_{i2}=\alpha_{i4}=1, \alpha_{i1}=0$, then $ \alpha_{i5}=\alpha_{i7}=2,\alpha_{i3}=\alpha_{i6}=\alpha_{i9}=\alpha_{i11}=\alpha_{i13}=\alpha_{i15}=1, \alpha_{in}=0$ for all $i \in \mathcal{I}, n \in \mathcal{T}\backslash\{2,3,4,5,6,7,9,11,13,15\}$ and $\beta_{i2}=\beta_{i3}=\beta_{i4}=\beta_{i5}=\beta_{in}=\beta_{i6}=\beta_{i7}=1, \beta_{in}=0$ for all $i \in \mathcal{I}, n \in \mathcal{T}\backslash\{2,3,4,5,6,7\}$  where $t^*_{i1}=3, t^*_{i2}=4$. The eighth case: $ \alpha_{i1}=\alpha_{i2}=\alpha_{i4}=1$, then $ \alpha_{i3}=\alpha_{i5}=\alpha_{i7}=2, \alpha_{i6}=\alpha_{i9}=\alpha_{i11}=\alpha_{i13}=\alpha_{i15}=1, \alpha_{in}=0$ for all $i \in \mathcal{I}, n \in \mathcal{T}\backslash\{1,2,3,4,6,5,7,9,11,13,15\}$ and $\beta_{i2}=\beta_{i3}=\beta_{i4}=\beta_{i5}=\beta_{i6}=\beta_{i7}=1, \beta_{in}=0$ for all $i \in \mathcal{I}, n \in \mathcal{T}$  where $t^*_{i1}=2, t^*_{i2}=3, t^*_{i3}=4$.

Conversely, given a feasible solution $\left(\mathbf{\alpha}, \mathbf{\beta}, \mathbf{t^*}\right)$ to the problem \eqref{NPHardModel}, we can construct a feasible solution $\{\mathcal{H}_l\}_{l \in \mathcal{L}}$ for the subset-sum problem by choosing  $\mathcal{H}_{1}=\{i: 2 \in \Gamma_i^* \}$, $\mathcal{H}_{2}=\{i: 3 \in \Gamma_i^* \}$  and  $\mathcal{H}_{4}=\{i: 4 \in \Gamma_i^* \}$. As $\left(\mathbf{\alpha}, \mathbf{\beta}, \mathbf{t^*}\right)$ is a feasible solution of problem \eqref{NPHardModel}, we satisfy the condition $\sum_{l\in \mathcal{L}}\sum_{i \in \mathcal{H}_l} w_{il}^{'}=W^{'}$. \Halmos
\endproof

\subsection{Proof of Proposition \ref{prop:binary}} \label{App:binary}

\textbf{Proposition 4} \textit{Let the state variables be binary with $x_{in} \in \{0, 1\}$ for all $i \in \mathcal{I}, \> n \in \mathcal{T}$. Then, the integrality conditions on the revision decision variables $r_{it}$ for all $i \in \mathcal{I}$ and $t \in \mathcal{P}$ can be relaxed within $(AMS_\mu)$.} 

\proof{Proof:}
For any given node pair $( m, n)$ from set $S_{t^{'}} \cap \mathcal{T}(l)$ where $l \in S_{t}, t^{'} = t + 1 , t \in \mathcal{P}, i \in \mathcal{I}$, the NACs  can be written as follows: 
\begin{subequations}\label{MMNLP}
\allowdisplaybreaks\begin{align}
&\hspace{0.5em} x_{i m} \geq x_{i n}-\left(r_{it^{'}}-r_{it} \right)  \\
&\hspace{0.5em} x_{i m} \leq x_{i n}+\left(r_{it^{'}}-r_{it} \right)
\end{align}
\end{subequations}
where $r_{it^{'}} -r_{it} \leq 1 $ and $r_{i1}= 0$. We distinguish two cases for the values of $x_{im}$ and $x_{in}$. \\
Case 1: If $x_{im} \neq x_{in}$ for at least one $(m,n)$ combination, then $x_{im}=0$ and $x_{in}=1$ or  $x_{im}=1$ and $x_{in}=0$. By construction, it is enforced that $r_{it'}=r_{it}+1$. \\
Case 2: If $x_{im}=x_{in}$ for all $(m,n)$, then  $x_{im} = x_{in} = 0$ or $x_{im} = x_{in} = 1$. 
Then, we can consider the following cases: i) $r_{it'}=r_{it}$, ii) $r_{it'} - r_{it} = 1$, and  iii) $0 < r_{it'} - r_{it} < 1$. The potential problematic case, which is case iii, implies that at most $\mu - 1$ revisions can be made. This case can only happen if $z(AMS_{\mu - 1}) = z(AMS_\mu)$, thus making $\mu - 1$ revisions optimal. In that case, it is still possible to construct an optimal solution with integer revision variables.

These cases show that one can obtain an optimal solution by only considering one by one increase for revision variables in consecutive stages, allowing elimination of the need to enforce integrality of revision variables. 
\Halmos
\endproof

\subsection{Proof of Proposition \ref{Prop:Elimination1}} \label{App:Elimination1}

\textbf{Proposition 5} \textit{Let $t \in \mathcal{P}$, $t^{'} \geq t$, and $l \in S_{t}$ be given. Let $K = |S_{t^{'}} \cap \mathcal{T}(l)|$. 
Without loss of generality, we can order the nodes in $S_{t^{'}} \cap \mathcal{T} (l)$ as $(n, n +1, \dots,n + K-1)$. 
Then, it is sufficient to include constraints \eqref{AdaptiveModel_d} and \eqref{AdaptiveModel_e} only for the consecutive node pairs (e.g., $(n,n+1),(n+1,n+2),\dots,(n+K-1,n)$) from the set $S_{t^{'}} \cap \mathcal{T}(l)$ as follows: 
\begin{align}
\begin{split}
 x_{in} \geq x_{i(n+1)} -\bar{x}_i\left(r_{it'}-r_{it} \right) \quad \quad n \in  S_{t^{'}} \cap \mathcal{T}(l),\quad l \in S_t,\quad t'\geq t,\quad t \in \mathcal{P},\quad i \in \mathcal{I} \\
\end{split}
\end{align}
We note that we should remove the ordered node $n + K$ from the set $S_{t^{'}} \cap \mathcal{T} (l)$ and add the constraint between $x_{i(n+K-1)}$ and $x_{in}$, but to enhance readability we slightly abuse notation in constraints \eqref{eq:Property1} by using node $n + K$ for node $n$.}

\proof{Proof:}
We prove this by conditioning on the value of $(r_{it'}-r_{it}) \in \mathbb{Z}_{+}$. 

\noindent 
\textit{Case 1:} $(r_{it'}-r_{it}) \geq 1$. Then both \eqref{AdaptiveModel_d} and \eqref{AdaptiveModel_e} are nonbinding, as are \eqref{eq:Property1}.

\noindent
\textit{Case 2:} $(r_{it'}-r_{it}) = 0$. The constraints in \eqref{eq:Property1} will reduce to $x_{in} = x_{i(n+1)} = \dots = x_{i(n+K-1)}$, where $K = |S_{t'} \cap \mathcal{T}(l)|$. Hence, for $m, n \in  S_{t^{'}} \cap \mathcal{T}(l)$ as defined in \eqref{AdaptiveModel_d} and \eqref{AdaptiveModel_e}, there exists unique $k, k' \in \{0, 1, \dots, K-1\}$ such that $x_{im} = x_{i(n + k)}$ and $x_{in} = x_{i(n + k')}$. 
Thus, constraints \eqref{eq:Property1} impose $x_{im} = x_{in}$ for all $m, n \in  S_{t^{'}} \cap \mathcal{T}(l)$ as defined in \eqref{AdaptiveModel_d} and \eqref{AdaptiveModel_e}. ~\Halmos 
\endproof

\subsection{Proof of Proposition \ref{Prop:Elimination2}} \label{App:Elimination2}

\textbf{Proposition 6} \textit{Let $t_{a(m,n)}$ be the stage of the last common ancestor of node pair $(m,n)$ from the set $S_t$ of stage $t$. It is sufficient to include NACs of $(m,n)$ only for $t, t_{a(m,n)}$ combination, and accordingly, constraints \eqref{AdaptiveModel_d} and \eqref{AdaptiveModel_e}  can be replaced by the following constraints:
\begin{subequations} \label{C:elimination2}
\begin{align}
& x_{im} \geq x_{in}-\bar{x}_i\left(r_{it}-r_{it_{a(m,n)}} \right) \quad \forall m, n \in S_{t},  t \in \mathcal{P}, i \in \mathcal{I} \\
& x_{im} \leq x_{in}+\bar{x}_i\left(r_{it}-r_{it_{a(m,n)}} \right) \quad \forall m, n \in S_{t},   t \in \mathcal{P},  i \in \mathcal{I}  
\end{align}\end{subequations}}

\proof{Proof:}
For a given $i$ and node pair $(m, n)$ from set $S_{t^{'}}$,  we have NACs \eqref{AdaptiveModel_d} and \eqref{AdaptiveModel_e} for some set of t, where $t\leq t'$, if $( m, n) \in \bigcup \limits_{l \in S_t} S_{t^{'}} \cap \mathcal{T}(l)$. Accordingly, we have NACs for $t$ and $t-v$ where $t=t_{a(m,n)}$ and $v \in \{1,\dots,t-1\}$  as follows:
\begin{subequations} \label{eq:Property2Proof}
\vspace{-4mm}
\allowdisplaybreaks\begin{align}
&\hspace{0.5em} x_{i  m} \geq x_{i n}-\bar{x}_i\left(r_{i  t'}-r_{it} \right)  \label{Prop2:a} \\
&\hspace{0.5em} x_{i m} \leq x_{i n}+\bar{x}_i\left(r_{i  t'}-r_{it} \right) \label{Prop2:b}  \\
&\hspace{0.5em} x_{i m} \geq x_{i n}-\bar{x}_i\left(r_{i  t'}-r_{i(t-v)} \right) \quad v \in \{1,\dots,t-1\} \label{Prop2:c}  \\
&\hspace{0.5em} x_{i m} \leq x_{i n}+\bar{x}_i\left(r_{i  t'}-r_{i(t-v)} \right) \quad v \in \{1,\dots,t-1\} \label{Prop2:d}
\end{align}
\end{subequations}
We can establish the relationship: $r_{i(t-v)} \leq r_{it} \leq r_{it'}$, where $t-v \leq t \leq t'$, based on constraint \eqref{AdaptiveModel_f}. This relationship implies that $(r_{it'} - r_{it}) \leq (r_{it'} - r_{i(t-v)})$. Since $(r_{it'}-r_{it})$ attains the lower value,  constraints \eqref{Prop2:a}--\eqref{Prop2:b} dominate constraints \eqref{Prop2:c}--\eqref{Prop2:d} for any $v \in \{1,\dots,t-1\}$, which shows that inclusion of only constraints \eqref{Prop2:a}--\eqref{Prop2:b} is sufficient for non-anticipativity. As a result, we conclude that for a node pair $(m,n)$ from set $S_{t^{'}}$, it is sufficient to include NACs only for $(t',t_{a(m,n)})$ combination. For notational brevity, we omit the prime symbol and consider $(t,t_{a(m,n)})$ combination in constraints \eqref{C:elimination2}.
\Halmos 
\endproof

\subsection{Proof of Proposition \ref{Prop:Elimination3}} \label{App:Elimination3}

\textbf{Proposition 7} \textit{Let $t\in \{1,\dots,T\}$ and $t'\geq t$. For given $\mu$, NACs between a particular $(t',t)$ combination become redundant, if $t'>t+(T-\mu)-1$. Then, constraints \eqref{AdaptiveModel_d} and \eqref{AdaptiveModel_e} reduce to:
\begin{subequations}
\allowdisplaybreaks\begin{align}
x_{i m} \geq x_{i n}-\bar{x}_i\left(r_{it'}-r_{it} \right) \ \forall m, n \in S_{t^{'}} \cap \mathcal{T}(l), l \in S_{t},  \min(t+ (T-\mu)-1, T) \geq t^{'} \geq t, t \in \mathcal{P}, i \in \mathcal{I} \label{eq:El31}\\
x_{i m} \leq x_{i n}+\bar{x}_i\left(r_{it'}-r_{it} \right) \  \forall m, n \in S_{t^{'}} \cap \mathcal{T}(l), l \in S_{t}, \min(t+ (T-\mu)-1, T) \geq t^{'} \geq t, t \in \mathcal{P}, i \in \mathcal{I} \label{eq:El32}
\end{align}
\end{subequations}}
 \vspace{-5mm}
\proof{Proof:}
From Lemma~\ref{Lem:binding}, at an optimal solution of $(AMS_\mu)$,  $r^*_{iT}=\mu$. 
Combining this observation with constraints \eqref{AdaptiveModel_f} - \eqref{AdaptiveModel_g} implies that $r_{iT}$ can reach $\mu$ only with a sequential increase of revision variables one by one. Consequently, if the difference between $t'$ and $t$ is greater than $T-\mu-1$, then the difference between $r_{it'}$ and $r_{it}$ is positive. Proof of Proposition \ref{Prop:Elimination1} shows redundancy of NACs if $r_{it'}-r_{it}>0$, therefore constraints \eqref{eq:El31} and \eqref{eq:El32} are sufficient for non-anticipativity. \Halmos
\endproof

\section{Scenario Tree Generation Algorithm}
\label{sec:AppendixScenarioTreeAlg}

\begin{algorithm}
\small
\caption{Scenario Tree Generation for Generation Expansion Planning Problem}
\label{ScenarioTreeGeneration}
\begin{algorithmic}[1]
\STATE Calculate the demand values at the root node as 
$\{d_{1k} : d_{1k}=w_k d_1 \}_{k \in \mathcal{K}}$ \\
\STATE Obtain capacity factors at root node as $ \{l_{i1k}\}_{i \in \mathcal{I}, k \in \mathcal{K}}$ from $\mathcal{N}(\mu^l_{i},\,\sigma^{l}_{i})$
\FORALL{$t \in \{2,\dots,T\}$}
\FORALL{$n \in S_t $}
\FORALL{$k \in \mathcal{K} $}
\STATE Generate growth factors $\lambda^d_{nk} \sim \mathcal{N}(\mu^d,\,\sigma^{d}) $
\STATE $d_{nk} \leftarrow (1 + \lambda^d_{nk}) d_{a(n),k} $ 
\FORALL{$i \in \mathcal{I} $}
\STATE $l_{ink} \sim \mathcal{N}(\mu^{l}_{i},\,\sigma^{l}_{i})$
\ENDFOR
\ENDFOR
\ENDFOR
\ENDFOR
\end{algorithmic}
\end{algorithm} 

\end{APPENDICES} 

\end{document}